\documentclass[a4paper, fleqn]{cas-dc}

\usepackage{amsthm}

\usepackage{subcaption} 
\usepackage{graphicx}
\usepackage{svg}
\usepackage{epstopdf}
\usepackage{caption}
\usepackage{float}
\usepackage{algorithm}
\usepackage{algorithmicx}
\usepackage{algpseudocode}
\usepackage{longtable}
\usepackage[figuresright]{rotating}
\usepackage{pdflscape}
\usepackage[hypcap = true, hypcapspace = 2cm]{caption}
\usepackage{makecell}
\usepackage{multicol}
\usepackage{enumitem}
\usepackage{bigstrut}
\usepackage{tabularx}
\usepackage{threeparttable} 

\usepackage{booktabs}
\usepackage{array, caption, threeparttable}
\newcolumntype{C}[1]{>{\centering\arraybackslash}p{#1}}

\usepackage{amsmath,bm}

\allowdisplaybreaks[4]
\makeatletter
\newenvironment{breakablealgorithm}
{
		\begin{center}
			\refstepcounter{algorithm}
			\hrule height.8pt depth0pt \kern2pt
			\renewcommand{\caption}[2][\relax]{
				{\raggedright\textbf{\ALG@name~\thealgorithm} ##2\par}%
				\ifx\relax##1\relax 
				\addcontentsline{loa}{algorithm}{\protect\numberline{\thealgorithm}##2}%
				\else 
				\addcontentsline{loa}{algorithm}{\protect\numberline{\thealgorithm}##1}%
				\fi
				\kern2pt\hrule\kern2pt
			}
		}{
		\kern2pt\hrule\relax
	\end{center}
}
\makeatother

\usepackage[section]{placeins} 

\newcommand{\mysubtitle}[2]{%
	\par\vspace{\baselineskip}        
	\noindent\hspace*{1.3em}(#1)\ \textit{#2}\par 
	\vspace{\baselineskip}            
}

\usepackage[many]{tcolorbox}
\usepackage{framed}
\usepackage{hyperref} 
\hypersetup{colorlinks=true,
            linkcolor=blue,   
            anchorcolor=black, 
            citecolor=blue,    
            urlcolor=blue      
            }

\usepackage[numbers, sort&compress]{natbib}
\usepackage[numbers]{natbib}

\usepackage{appendix}
\renewcommand\appendix{\par
    \setcounter{section}{0}
    \setcounter{subsection}{0}
    \setlength\LTleft{0pt}
    \setlength\LTright{0pt}
    \captionsetup[table]{margin=-0.03cm}
    \gdef\thesection{Appendix \Alph{section}}}

\def\tsc#1{\csdef{#1}{\textsc{\lowercase{#1}}\xspace}}
\tsc{WGM}
\tsc{QE}

\begin{document}
\let\WriteBookmarks\relax
\renewcommand{\floatpagefraction}{0.8}
\renewcommand{\textfraction}{0.05}
\let\printorcid\relax 
\shorttitle{}



\title[mode = title]{The multi-objective portfolio model for oil and gas exploration drilling projects selection and its operator-enhanced NSGA-II based solution}  




\author[1,2,3]{Chao Min}[style=chinese]
\author[1, 2]{Junyi Cui}[style=chinese]
\author[3]{Stanis{\l}aw Mig{\'o}rski}[style=chinese]
\cormark[1]
\author[1]{Yonglan Xie}[style=chinese]
\author[1]{Qingxia Zhang}[style=chinese]
\author[4]{Jun Peng}[style=chinese]

\address[1]{School of Sciences, Southwest Petroleum University, Chengdu 610500, China}
\address[2]{Institute for Artificial Intelligence, Southwest Petroleum University, Chengdu 610500, China}
\address[3]{Jagiellonian University in Krakow, Faculty of Mathematics and Computer Science, 30348 Krakow, Poland}
\address[4]{School of Geoscience and Technology, Southwest Petroleum Uniersity, Chengdu, Sichuan 610500, China}

\cortext[1]{Jagiellonian University in Krakow, Faculty of Mathematics and Computer Science, 30348 Krakow, Poland. E-mail address: stanislaw.migorski@uj.edu.pl.} 

\begin{abstract}
Drilling investment is pivotal to operational planning in oil and gas (O\&G) exploration. Conventional deployment relies heavily on fragmented expert assessments of geological and economic factors, with limited integration ability of information. As the tool of portfolio show strong potential for mitigating uncertainty and selecting superior drilling plans, this study develops a multi-objective mean–variance portfolio model that accounts for geological-parameter uncertainty, enabling an effective risk–return trade-off and optimal selection. First, the probabilistic distribution of geological-parameters for prospect-list projects is obtained through expert-elicited priors. And considering the selection of the drilling projects as a portfolio, an optimization model is formulated jointly to describe the return and risk of short-term plan, under different constraints. Second, an improved OE-NSGA-II algorithm is proposed specifically for this model, in which (i) a directional crossover operator is designed to embed improving directions in objective space—derived from dominance and objective differences—into recombination, and (ii) a structure-aware mutation operator is designed to prioritize high-utility bit flips via probabilistic sampling with feasibility repair, thus improving the search ability for superior Pareto solutions. Finally, using the case of  2023 exploration drilling deployment for verification, and then apply the validated method to the 2024 deployment to support decision-making. The results indicate that the proposed approach offers a reusable solution for drilling portfolio optimization in O\&G exploration.
\end{abstract}


\maketitle 

\section{Introduction}
O\&G exploration is characterized by capital intensity, high risk, and long cycles. Decision-makers must rely on subjective geological assessments and economic risk evaluations to compile drilling investment prospect inventory (i.e., candidate projects for drilling) \cite{bib1}. In China, oil companies typically allocate capital across geophysical prospecting, traps, and appraisal. Within stage-based investment strategies, preliminary prospected wells and appraisal wells receive prioritized funding. This study focuses on the drilling investments for trap and appraisal projects within a basin’s licensed area and employs optimization techniques to generate actionable shortlists of preferred drilling projects, providing quantitative support for resource allocation and risk control.

Since the 1980s, Markowitz’s portfolio theory \cite{bib2} has been introduced into project investment in O\&G industry to investigate its applicability in the energy sector—particularly in exploration and development—and has played a significant role. Thereafter, many scholars lay attention to the balance between expected return and its variance/semi-variance in exploration–development portfolios to characterize the trade-off between return and risk. For example, Hightower et al. \cite{bib3} applied an improved portfolio model to the selection of exploration drilling plans in petroleum companies, offering a new perspective for drilling decision-making. Guo et al. \cite{bib4} adopts net present value (NPV) as the value metric and semi-variance as the risk measure, demonstrating the model’s effectiveness and practicality in real applications. Martijn et al. \cite{bib5} proposes a portfolio optimization method for O\&G projects under budget and capacity constraints, grounded in quadratic programming and preference theory. Tang et al. \cite{bib6} investigates key elements of portfolio design by incorporating covariances among prospects (projects), thereby constructing a more robust optimization framework. Jorge et al. \cite{bib7} examines mean–semivariance portfolio selection for projects, using the semivariance of portfolio NPV as risk and the expected portfolio NPV as return. Yan et al. \cite{bib8} models cash flows with normal uncertain variables and proposes maximizing expected return while minimizing the sine cross-entropy between realized and prior returns. an et al. \cite{bib9} develops a nonlinear multi-objective mixed-integer programming model for O\&G portfolios that accommodates nonlinear equations and integer constraints. 
\begin{tcolorbox}[%
	title=Nomenclature,
	float*=t,
	breakable,             
	colframe=black,        
	colback=white,         
	boxrule=1pt,
	left=1mm, right=1mm,   
	before skip=1ex,       
	after skip=1ex,
	width=\textwidth]
	\begin{multicols}{2}
		\noindent \textit{Abbreviations}
		\begin{description}
			\item[O\&G]           Oil and gas
			\item[MV]           Mean-variance
			\item[GPoS]         Probability of geological success
			\item[EPoS]         Probability of economic success
			\item[MEFS]          Minimum economic field size
			\item[NPV]           Net-present value
			\item[EMV]           Expected monetary value
		\end{description}
		
		\vspace{1ex}
		
		\noindent\textit{Subscripts}
		\begin{description}
			\item[\textit{t}]    Year
			\item[\textit{oi}]   Oil
			\item[\textit{ga}]    Gas
			\item[\textit{tra}]  Trap project
			\item[\textit{app}]  Appraisal project
			\item[\textit{i}]    Trap projects resources
			\item[\textit{j}]    Appraisal projects resources
			\item[\textit{k}]    $ k \in \left\{ {0,\;1} \right\} $
			\item[\textit{lb}]    Lower bound value
			\item[\textit{ub}]    Upper bound value
		\end{description}
		
		\vspace{1ex}
		
		\noindent\textit{Parameters}
		\begin{description}
			\item[$ \textit{PoS}_{i}^{g} $]  GPoS of trap project $ i $
			\item[$ \textit{PoS}_{j}^{e} $]  EPoS of appraisal project $ j $
			\item[$ \textit{Pred}^{oi} $]  The predicted oil reserves
			\item[$ \textit{Pred}^{ga} $]  The predicted gas reserves
			\item[$ \textit{Cont}^{oi} $]  The controlled oil reserves
			\item[$ \textit{Cont}^{ga} $]  The controlled gas reserves
			\item[$ \textit{Prov}^{oi} $]  The proved oil reserves
			\item[$ \textit{Prov}^{ga} $]  The proved gas reserves
			\item[$ \textit{Cost}_{i}^{tra} $]  The drilling costs for trap project $ i $
			\item[$ \textit{Cost}_{j}^{app} $]  The drilling costs for appraisal project $ j $
			\item[$ \textit{Cost}_{ub}^{tra} $] The investment cost of trap projects
			\item[$ \textit{Cost}_{ub}^{app} $] The investment cost of appraisal projects
			\item[$ \textit{Drill}^{tra} $] The number of wells drilled in the trap project
			\item[$ \textit{Drill}^{app} $] The number of wells drilled in the appraisal project
			\item[$ \textit{Tot}_{wells} $]      Total number of wells drilled
			\item[$ \textit{Thre}_{well} $]      Threshold for the number of low probability of success for a single well
			\item[$ \textit{pred}_{lb }^{oi} $]    Lower bound of predicted oil reserves
			\item[$ \textit{pred}_{lb }^{ga} $]    Lower bound of prospective gad reserves
			\item[$ \textit{Cont}_{lb }^{oi} $]    Lower bound of controlled oil reserves
			\item[$ \textit{Cont}_{lb }^{ga} $]    Lower bound of controlled gad reserves
			\item[$ \textit{Prov}_{lb }^{oi} $]    Lower bound of proved oil reserves
			\item[$ \textit{Prov}_{lb }^{ga} $]    Lower bound of proved gad reserves
			\item[$ {L}_{ub} $]         Upper limit on the number of projects whose geological or economic probability of success meets or exceeds the threshold
			\item[\textit{$ N_{tra}^{\left( r \right)} $}]    Threshold for the number of selected trap projects per region
			
			\item[\textit{$ N_{app}^{\left( r \right)} $}]   Threshold for the number of selected appraisal projects per region
		\end{description}
		
		\vspace{1ex}
		
		\noindent\textit{Decision variables}
		\begin{description}
			\item[$ x_{ik} $]   Whether the $ i^{th}$ trap project is selected
			\item[$ x_{jk} $]   Whether the $ j^{th}$ appraisal project is selected
		\end{description}
		
	\end{multicols}
\end{tcolorbox}

Existing studies have focused on refining model forms or modeling specific variables (particularly economic factors), but have rarely provided an integrated treatment of the geological, resource, and economic dimensions to capture the inter-variable uncertainty.

Some scholars argue that the mean–variance model remains one of the simplest and most effective tools for asset selection \cite{bib10,bib11}, it is important to recognize that, relative to the financial sector, applying modern portfolio theory in the oil and gas industry differs substantially in terms of the types of uncertainties, risk indicators, market characteristics, time periods, and the effect of budget constraints \cite{bib12,bib13}. Motivated by these gaps, this study develops a multi-objective optimization model that jointly incorporates geological, resource, and economic uncertainties, with the aim of delivering effective and reliable decision support for drilling deployment in O\&G exploration.

The proposed model is rooted in portfolio theory and can be formulated as a constrained multi-objective optimization problem. Using augmented Lagrangian methods to handle such equality constraints is a viable strategy. For example, Neven et al. \cite{bib14} formulate dynamic reconfiguration of distribution networks (DRDN) as a mixed-integer linear programming (MILP) problem and solve it via Lagrangian relaxation; Zhao et al. \cite{bib15} developed a parallel generalized Lagrangian–Newton solver for a class of PDE-constrained optimization problems with inequality constraints, and Wang et al. \cite{bib16} studied non-negative cardinality constraints and proposed a Lagrangian–Newton algorithm (LNA) for inequality-constrained problems.

However, in the context of drilling portfolio optimization, penalty terms will inevitably distort the objective values to some extent, which can interfere with non-dominated sorting and crowding-distance calculations, hindering the precise characterization of the return–risk trade-off. In addition, classical optimization techniques typically suffer from high computational cost and unsatisfactory convergence when tackling such NP-hard problems. 

In recent years, meta-heuristic algorithms have attracted widespread attention for this class of problems due to their ease of implementation, algorithmic flexibility, and strong ability to escape local optima \cite{bib17}. For example, Ku{\c{s}}o{\u{g}}lu et al. \cite{bib18} proposes a multi-objective Harris Hawks Optimizer (HHO); Zhao et al. \cite{bib19} develops a multi-objective artificial hummingbird algorithm (MOAHA) with dynamically pruned crowding distance; Sidi et al. \cite{bib20} introduces the multi-objective flower pollination algorithm (MOPFA); and Zhao et al. \cite{bib21} presents a multi-objective firefly algorithm with adaptive region division (MOFA-ARD). In the latest portfolio-optimization studies, Song et al. \cite{bib22} proposes a dual-time, dual-population MOEA (DTDP-EAMO) to reduce risk and increase return; Sa{\"i}b et al. \cite{bib23} develops an improved squirrel search algorithm (ISSA) to balance return and risk; Minh et al. \cite{bib24} couples multi-objective PSO (MOPSO) with a shrinkage covariance estimator, yielding a MOPSO–Shrinkage hybrid for volatile markets; Li et al. \cite{bib25} applies a multi-strategy modified sparrow search algorithm (MSMSSA) to portfolio selection with uncertain tail VaR; and Zhan et al. \cite{bib26} designs an enhanced dung beetle optimizer (EDBO) for multi-period portfolio problems under multi-factor uncertainty, demonstrating both effectiveness and practicality.

In general, NSGA-II \cite{bib27}—with its efficient non-do\-mi\-na\-ted sorting and crowding distance—has been widely adopted to obtain high-quality Pareto fronts in multi-objecti\-ve optimization. Nevertheless, for domain-specific problems with special structure, the classical NSGA-II may be inadequate. Tailored to our drilling portfolio model, we enhance NSGA-II’s crossover and mutation by proposing a parent-decision-bit–guided directional crossover (DC) and a structure-aware mutation with bit-level collisions (SAM), thereby improving the search efficiency and solution quality of OE-NSGA-II.

The main contributions are as follows:
\begin{itemize}
	\item The Iman–Conover method is employed to characterize the dependence among sub-factors of the probability of geological success (GPoS) and inject this dependence into the joint samples, thereby enhancing the reliability of the tail quantiles of the GPoS.
	\item Multiple sources of geological uncertainty are integrated into a mean–variance framework, yielding a multi-objective portfolio model for O\&G exploration drilling.
	\item An improved OE-NSGA-II, combining DC and SAM, is developed to preserve diversity while directionally enhancing the search for high-quality non-dominated solutions, thereby improving the return–risk trade-off.
\end{itemize}

The structure of this paper is organized: Section 2 presents the exploration background and defines the key parameters; Section 3 analyzes the parameter models and introduces a drilling portfolio optimization model under geological-parameter uncertainty; Section 4 refines the mod\-el and proposes the OE-NSGA-II algorithm; Section 5 applies the model and algorithm to case studies in O\&G exploration to validate their effectiveness; Section 6 concludes the paper and outlines potential avenues for improvement.

\section{Description and definitions}

\subsection{Background}
When planning short-term exploration, oil companies must select preliminary prospected wells and appraisal wells from different blocks (trap projects and appraisal projects) for investment. In general, trap projects are characterized by incomplete exploration data (geological, resource, and economic information) and require preliminary prospected wells drilling to determine block value. Appraisal projects, by contrast, are supported by relatively rich data and well-defined prospects, and can be further divided into drilling appraisal wells (which require drilling to assess commercial value) and reserve-providing appraisal wells (which can be evaluated without drilling). At present, investment decisions are largely driven by expert experience and qualitative judgment, with limited support from systematic deployment models or economic-theory-based frameworks. As a result, under uncertain geological success probabilities and concurrent constraints on reserves and capital, it is difficult to quantitatively balance “input–output” and achieve optimal allocation.

\begin{figure*}[htbp]
	\centering
	\includegraphics[width=\linewidth]{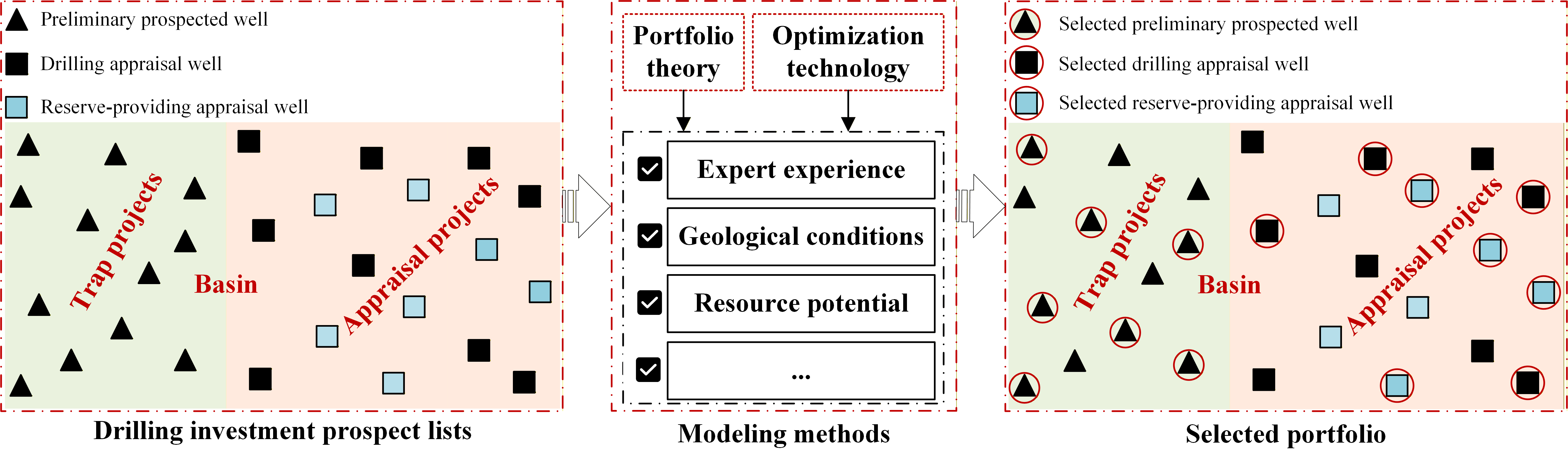} 
	\caption{Optimal process for O\&G exploration drilling.}
	\label{fig1}
\end{figure*}

To address this issue, we recast drilling deployment as a portfolio selection problem: each well (preliminary prospected well and appraisal well) is treated as an asset, with its expected NPV representing “return” and the variance of expected return characterizing the portfolio risk. On this basis, we formulate a multi-objective optimization model that supports quantitative decisions on “which wells to drill and at what investment level.” A brief overview of the research background is provided in Fig. \ref{fig1}.

\subsection{Definitions}
Terminology in the petroleum industry is often ambiguous or varies between companies \cite{bib1}. This study unifies the key parameters as follows:

\textit{Drilling success}: the discovery of an oil well following the drilling of a preliminary prospected well or appraisal well.

\textit{Probability of geological success} (GPoS): the probability of encountering freely flowing, continuously producible petroleum from a penetrated formation under cautious operational conditions \cite{bib28}.

\textit{Geological success rate}: the proportion of drilled wells that yield oil. Ideally, for a large number of exploration wells, the mean geological success probability equals the actual geological success rate \cite{bib29}.

\textit{Success-case volume}: the expected volume of recoverable petroleum fluids if a well encounters movable petroleum with sustained flow. Most companies estimate the distribution of recoverable volumes via Monte Carlo simulation by multiplying distributions of total rock volume, net-to-gross ratio, porosity, hydrocarbon saturation, formation volume factor, and recovery factor \cite{bib1, bib30, bib31}. Here, P90, P50, and P10 denote reserves at the 90\%, 50\%, and 10\% confidence levels, respectively; Pmean represents the arithmetic mean, characterizing the center of the probability distribution \cite{bib1}.

\textit{Economic probability of success} (EPoS): the pre‐drill probability of discovering a minimum economic field size (MEFS) and achieving profitability based on technical and economic evaluations \cite{bib1, bib32}.

\textit{Minimum economic field size} (MEFS): the minimum volume of recoverable O\&G necessary to make a project an economic success \cite{bib32}.

\textit{Net present value} (NPV): the difference between the present value of cash inflows (e.g., revenues from oil sales) and the present value of cash outflows (e.g., expenditures on field facilities and operations) \cite{bib1}.

\textit{Expected monetary value} (EMV): the average monetary outcome estimated based on NPV, EPoS, drilling investment, and operating costs \cite{bib33}.

\section{Methodology}
\begin{figure*}[htbp]
	\centering
	\includegraphics[width=0.9\linewidth]{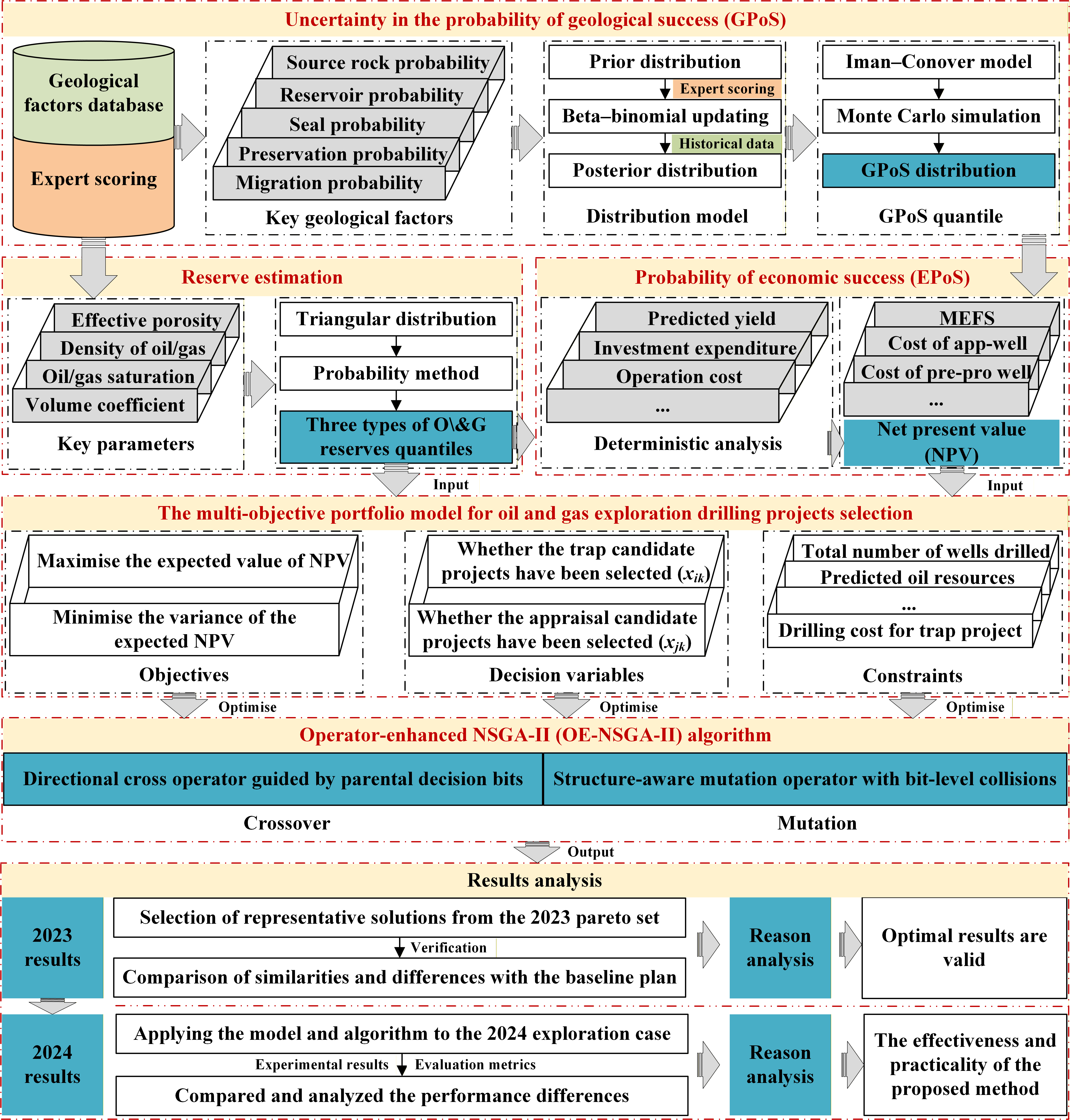} 
	\caption{Overall flowchart of the proposed model and algorithm.}
	\label{fig2}
\end{figure*}

This section presents all the necessary input parameters in sequence: the uncertainty analysis and quantification of the geological success probability, reserve estimation, the probability of economic success (EPoS), and the construction of multi-objective portfolio model for O\&G exploration drilling projects selection. The technical roadmap of the model and solution algorithm is shown in Fig. \ref{fig2}.

\subsection{Geological uncertainty factors}
The drilling portfolio optimization problem can be formulated in terms of an objective function
\begin{flalign}\label{1}
	& f:\;\mathscr{X} \times {\mathbb{R}^m} \to \mathbb{R},\;\left( {x,y} \right) \mapsto f\left( {x,y} \right), &
\end{flalign}
where $ x \in \mathscr{X} \subset {\mathbb{R}^n} $ denotes the decision vector and $ y \in {\mathbb{R}^m} $ is a realization of the random input vector $ Y $. The uncertainties considered in this study include the $\textit{PoS}^{g} $, $ \textit{PoS}^{g} $, three types of O\&G reserves (predicted oil reserves (Pre-OR), predicted gas reserves (Pre-GR), controlled oil reserves (COR), controlled gas reserves (CGR), proved oil reserves (Pro-OR), proved gas reserves (Pro-GR)), and NPV. Since the setting of the above parameters relies heavily on expert experience, this paper employs a probabilistic approach (Monte Carlo simulation) to quantify the associated uncertainties.

\subsection{Quantification of GPoS}
\mysubtitle{1}{Probability distribution of sub-factors for GPoS}

Determining GPoS imposes stringent requirements on the evaluation of subclass probabilities. It necessitates gathering historical exploration data—such as source rocks $ p_{1} $, reservoir $ p_2 $, preservation $ p_3 $, seal $ p_4 $, migration $ p_5 $—to quantify this critical parameter, with expert elicitation serving as the primary input for constructing the probability of geological success database \cite{bib34}. Let $ D = \left\{ {1,\;2,\; \ldots ,\;d} \right\} $ denote key geologic factors $ p_i $. The factor-level success probabilities are $ {p_i} \in \left( {0,1} \right) $ and the overall GPoS is $ \textit{PoS}^{g} = \prod\limits_{i = 1}^5 {{P_i}} $ \cite{bib35, bib36}. The values of these key geological parameters are typically derived from multiple experts’ judgments. Aggregating their assessments yields the minimum $ a $, mode (most likely value) $ b $, and maximum $ c $.  To represent the uncertainty, we model each parameter with a triangular distribution supported on $ \left[ {a,c} \right] $ with mode $ b $. Let $ p \in \left( {0,1} \right) $ be a uniform random variate; the inverse cumulative distribution function is:
\begin{flalign}\label{2}
	& f\left( p \right) = \left\{ \begin{array}{l}
		a + \sqrt {p\left( {c - a} \right)\left( {b - a} \right)} ,\;\;\;\;\;\;\;\;0 < p < \frac{{b - a}}{{c - a}} \\ 
		c - \sqrt {\left( {1 - p} \right)\left( {c - a} \right)\left( {c - a} \right)} ,\;\frac{{b - a}}{{c - a}} < p < 1 \\ 
	\end{array} \right. &
\end{flalign}
where the mean of this distribution is $ \mu  = \frac{{a + 4c + b}}{6} $ \cite{bib37,bib38}. 

These parameters are interdependent random variables, then the method in \cite{bib39} can be used to model inter-variable dependence, and \cite{bib40} has applied this idea to uncertainty sampling for the best estimate "standard" project (PWR12-BE) to estimate schedule, overnight capital cost (OCC), and total capital investment cost (TCIC). In contrast, to assess uncertainty for future observations, a Beta–Binomial framework is adopted to coherently fuse expert knowledge with sparse historical data. For each factor $ i $, a prior confidence (concentration) $ k_i $ is set to reflect data availability at the current exploration stage, and a Beta prior is placed on the success probability: $ {p_i} \sim Beta\left( {\alpha _i^{\left( 0 \right)},\beta _i^{\left( 0 \right)}} \right) $. Using a PERT aggregation of experts’ three-point assessments to obtain the prior mean $ {\mu _{i,\;\textit{PERT}}} $, the concentration parameters are:
\begin{flalign}\label{3}
	& \begin{array}{l}
		\alpha _i^{\left( 0 \right)} = {\mu _{i,\:\text{PERT}}}\left( {{k_i} - 2} \right) + 1, \; {k_i} > 0,\\ 
		\beta _i^{\left( 0 \right)} = \left( {1 - {\mu _{i,\:\text{PERT}}}} \right)\left( {{k_i} - 2} \right) + 1,\;{k_i} > 0. \\ 
	\end{array}
	&
\end{flalign}

Given historical data, let $ s_i $ and $ f_i $ denote the observed counts of successes and failures for factor $ i $. By conjugacy, the posterior is:
\begin{flalign}\label{4}
	& \begin{array}{l}
		p_i \mid \text{data} \sim \operatorname{Beta}\!\big(\alpha_i^{(1)},\, \beta_i^{(1)}\big),\: \\ 
		\alpha _i^{\left( 1 \right)} = \alpha _i^{\left( 0 \right)} + {s_i},\: \\ 
		\beta _i^{\left( 1 \right)} = \beta _i^{\left( 0 \right)} + {f_i}, \\ 
	\end{array}
	&
\end{flalign}
with posterior mean and marginal posterior given by:
\begin{flalign}\label{5}
	& \begin{array}{l}
		\mathbb{E}[\,p_i \mid \text{data}\,]
	= \frac{\alpha_i^{(1)}}{\alpha_i^{(1)}+\beta_i^{(1)}},\: \\
	p_i \mid \text{data} \sim \operatorname{Beta}\!\big(\alpha_i^{(1)},\,\beta_i^{(1)}\big). 
     \end{array}&
\end{flalign}

Using historical data, the Spearman rank-correlation matrix $ {R_S} $ is computed firstly. Because the empirical correlation matrix may be non–positive semidefinite (PSD), it is projected onto the nearest PSD correlation matrix as follows: perform the eigendecomposition $ {R_S} = Q\Lambda {Q^{\rm T}} $, apply eigenvalue clipping with $ \varepsilon  > 0 $: $ {\Lambda _\varepsilon } = diag\left( {\max \left\{ {{\lambda _i},\varepsilon I} \right\}} \right) $, and form $ {B} = Q\Lambda_{\varepsilon} {Q^{\rm T}} $. Let $ D = \textit{diag}{\left( B \right)^{ - \frac{1}{2}}} $ (where $ \textit{diag}\left( B \right) $ extracts the diagonal and forms a diagonal matrix); the normalized matrix $ {\hat R_S} = DBD $ satisfies $ {\hat R_S} \succeq O $ and $ \textit{diag}\left( {{{\hat R}_S}} \right) = 1 $. With $ \varepsilon  > 0 $, $ {\hat R}_S $ can also be guaranteed to be strictly positive definite.

\mysubtitle{2}{Monte Carlo and Iman-Conover methods}
For each factor $ i $, independently draw $ N $ samples from its posterior distribution and stack them into:
\begin{flalign}\label{6}
	& {X^{\textit{ind}}} = \left[ {{p_1}\;{p_2}\;{p_3}\;{p_4}\;{p_5}} \right] \in {\left( {0,1} \right)^{N \times d}}, &
\end{flalign}
apply the Iman–Conover procedure (rank matching with correlation projection) to $ X^{\textit{ind}} $ to obtain $ X^{\textit{corr}} $ whose (Spearman) correlation matches the target $ {\hat R_S} $.

For the $ t $-th Monte Carlo samples, let the correlated probability vector be $ {p^{\left( t \right)}} = \left( {p_1^{\left( t \right)},\;p_2^{\left( t \right)},\;p_3^{\left( t \right)},\;p_4^{\left( t \right)},\;p_5^{\left( t \right)}} \right) $. Define the combined probability $ p_g^{\left( t \right)} = g\left( {{p^{\left( t \right)}}} \right) $. This yields $ \left\{ {p_g^{\left( t \right)}} \right\}_{t = 1}^N $ over $ N $ simulations.

The sample mean and sample standard deviation of the combined probability are:
\begin{flalign}\label{7}
	& \begin{array}{l}
		\hat \mu  = \frac{1}{N}\sum\limits_{n = 1}^N {{P_g}^{\left( n \right)}} ,\;\\
		\hat \sigma  = \sqrt {\frac{1}{{N - 1}}\sum\limits_{n = 1}^N {{{\left( {{P_g}^{\left( n \right)} - \hat \mu } \right)}^2}} }.  
	\end{array}&
\end{flalign}

In summary, we decompose the sub-probabilities and combine expert judgment with historical data, employing Monte Carlo–IC sampling to simulate each random sub-probability and explicitly inject inter-factor dependence into the joint samples of $ \textit{PoS}^{g} $. This yields a more faithful representation of overall uncertainty, enhances the reliability of the tail quantiles of $ \textit{PoS}^{g} $, and provides a sound input basis for quantifying return and risk in the subsequent portfolio model. The estimation of three types of O\&G reserves is further developed below within this framework.

\subsection{Reserve estimation}
Uncertainty in reserves estimation mainly arises from the evaluation of geological data, which in turn depends heavily on expert judgment and field experience. Consequently, analogy method \cite{bib41}, volumetric method \cite{bib42}, and probability method \cite{bib43} have become standard techniques for estimating oil and gas reserves. This study focuses on carbonate fractured–cave reservoirs. Geng et al. \cite{bib43} constructed probability distribution models for carbonate reservoirs and proposed the PRUC model, which offers valuable insight for reserves estimation. Building on experts’ three-point assessments of effective porosity and oil saturation, we use the methodology described in the previous section to define triangular distribution types, and then apply Monte Carlo simulation to obtain distributed reserve estimates for candidate projects across different exploration stages and regions. The main methods for reserves estimation are briefly summarized below.{\tiny }
\begin{flalign}\label{8}
	& {G_{oi}} = 100\Phi {S_{oi}}{\rho _{oi}}/{\beta _{oi}}, &
\end{flalign}
where $ G_{oi} $ RPUC of crude oil, $ {10^4}t/k{m^2} $; $ \Phi $ effective porosity; $ \rho _{oi} $ the density of oil on the ground, $ g/c{m^3} $; $ S_{oi} $ oil saturation; $ \beta _{oi} $ volume coefficient of crude oil.
\begin{flalign}\label{9}
	&{G_{ga}} = 0.01\Phi {S_{ga}}{\rho _{ga}}/{\beta _{ga}},&
\end{flalign}
where $ G_{ga} $ denotes recoverable gas reserves $ ({10^8}m^3/k{m^2}) $;  $ \Phi $ effective porosity; $ \rho _{ga} $ the density of gas on the ground, $ g/c{m^3} $; $ S_{ga} $ gas saturation; $ \beta _{ga} $ volume coefficient of gas.

\subsection{Quantification of EPoS}
EPoS is defined as the product of GPoS and the probability of exceeding the MEFS, representing the likelihood that the discovered cumulative resources exceed the MEFS threshold. It is calculated as: 
\begin{flalign}\label{10}
	& \textit{PoS}^{e} = \textit{PoS}^{g} \times {P_{\textit{MEFS}}}. &
\end{flalign}

In this study, NPV is used primarily to discount cash flows associated with drilling, and is calculated as follows:
\begin{flalign}\label{11}
	&\textit{NPV} = \sum\limits_{t = 0}^n {\left( {{\textit{CI}} - {\textit{CO}}} \right)} {\left( {1 + r} \right)^{ - t}},&
\end{flalign}
where, $ \textit{CI} $ and $ \textit{CO} $ denote cash inflows and cash outflows (\textit{$ 10^4 $ } CNY), respectively. $ r $ is the discount rate. $ t $ is the evaluation period in years. Therefore, accounting for investment success and failure, the EMV is calculated as:
\begin{flalign}\label{12}
	&\textit{EMV} = \left( {\textit{NPV} \times \textit{PoS}^e} \right) - \left[ {\textit{Costs} \times \left( {1 - \textit{PoS}^e} \right)} \right].&
\end{flalign}

\subsection{The multi-objective portfolio model for O\&G exploration drilling projects selection}
This subsection details the multi-objective portfolio model for O\&G exploration drilling projects selection, including the model assumptions, decision variables, objective functions, and constraints.

\subsubsection{Assumptions}
This model involves multiple influencing factors. For modeling and solution purposes, the following assumptions are made:
\begin{itemize}
	\item All candidate projects have undergone preliminary assessment and screening, with technical feasibility and economic viability provisionally confirmed; project returns can be estimated using historical data, industry expertise, and simulation.
	\item Exploration funding is assumed to derive solely from internal capital, with no consideration of fundraising costs or financing arrangements.
	\item If an exploration target is confirmed to have economic value, its development is assumed to proceed unimpeded by external factors and to generate sustained positive cash flows post-development.
	\item Decision-makers are assumed to be rational agents seeking to balance maximization of investment returns and minimization of risk.
	\item Investment proportions for the two project types can be allocated between 0\% and 100\%, with the sum constrained to 100\%.
\end{itemize}

\subsubsection{Decision variables}
Let $ A $ denote the set of all candidate wells, including both trap and appraisal projects. The subset $ B = \left\{ {{x_{1k}},\;{x_{2k}}, \ldots ,{x_{mk}}} \right\} $ represents the preliminary prospected wells in trap projects, and the subset $ C = \left\{ {{x_{m + 1,k}},\:{x_{m + 2,k}},} \right. $ $ \ldots , \left. {{x_{nk}}} \right\} $ represents the appraisal wells in appraisal projects. For each candidate well, we define a binary decision variable $ x_{ik} \in \left\{ {0,\;1} \right\} $, where $ x_{ik} =1 $ indicates that the $ i $-th preliminary prospected wells or appraisal wells is selected. In practical drilling deployment, a subset of preliminary prospected wells and appraisal wells must be designated as mandatory to support risk investment in new blocks, and the corresponding decision variables are therefore fixed to 1. For notational convenience, all decision variables are stacked into a single vector $ z $:
\begin{flalign}\label{13}
	&z = \left[ {{x_{1k}} \ldots \;{x_{ik}} \ldots {x_{mk}}\;{x_{m + 1,k}} \ldots \;{x_{jk}} \ldots \;{x_{nk}}} \right].&
\end{flalign}

\subsubsection{Objective function}
In the mean–variance framework, this study fully integrates geological, economic, and reservoir engineering factors with the objective of selecting the most prospective drilling portfolio by maximizing the EMV while minimizing risk (variance). Given that trap projects are in the reconnaissance and preliminary exploration stages with uncertain geological data, geological success probability is used to quantify their success likelihood. Thus, the \textit{EMV} for trap projects is expressed as:
\begin{flalign}\label{14}
	&\begin{array}{l}
		\max : \\
		\textit{TR} = \sum\limits_{{ik} \in B} {\left[ {\left( {\textit{NPV}{_{ik}} \times \textit{PoS}_{ik}^g} \right) - \textit{Costs}_{ik}^{tra}} \right]},
	\end{array}&
\end{flalign}

The first term on the right-hand side of the equation represents the expected monetary value of each drilling project in the portfolio, where $ \textit{Costs}_{ik} $ denotes the drilling cost. Appraisal projects are at the appraisal stage with clear geological conditions and rely primarily on economic evaluation; thus, their \textit{EMV} can be expressed as:
\begin{flalign}\label{15}
	&\begin{array}{l}
		\max : \\ 
		\textit{AP} = \sum\limits_{jk \in C} {\left[ {\left( {\textit{NPV}{_{jk}} \times \textit{PoS}_{jk}^e} \right) - \textit{NPV}{_{jk}} \times \left( {1 - \textit{PoS}_{jk}^e} \right)} \right]},  \\ 
	\end{array}&
\end{flalign}
here, the first term on the right-hand side denotes the expected monetary value of each drilling project in the portfolio based on economic evaluation, and the second term represents the expected loss from failed cases. Therefore, the $ \textit{EMV} = \textit{TR} + \textit{AP} $ for O\&G exploration drilling portfolio optimization can be expressed as:
\begin{flalign}\label{16}
	&\begin{array}{l}
		\max :\:\: \\
		\textit{EMV} = \sum\limits_{{ik} \in B} {\left[ {\left( {\textit{NPV}{_{ik}} \times \textit{PoS}_{ik}^g} \right) - \textit{Costs}{_{ik}}} \right]}  \\
		\;\;\;\;+ \sum\limits_{{jk} \in C} {\left[ {\left( {\textit{NPV}{_{jk}} \times \textit{PoS}_{jk}^e} \right) - \textit{NPV}{_{jk}} \times \left( {1 - \textit{PoS}_{jk}^e} \right)} \right]}.
	\end{array}&
\end{flalign}

Similarly, the mathematical expression for risk (variance) is given by:
\begin{flalign}\label{17}
	&\begin{array}{l}
		\min :\:\: \\ 
		Risk = \left[ {{{\sum\limits_{ik \in B} {\left( {\textit{NPV}_{ik}\cdot \textit{PoS}_{ik}^g - \mu \left( {{x_{ik}},\:{x_{jk}}} \right)} \right)} }^2}\cdot{x_{ik}}} \right. \\ 
		\;\;\;\;\;\;\;\;{\left. { + {{\sum\limits_{jk \in C} {\left( {\textit{NPV}_{jk}\cdot \textit{PoS}_{jk}^e - \mu \left( {{x_{ik}},\:{x_{jk}}} \right)} \right)} }^2}\cdot{x_{jk}}} \right]^{\frac{1}{2}}}, \\ 
	\end{array}&
\end{flalign}
where $ \mu \left( {{x_{ik}},\;{x_{jk}}} \right) = \frac{{\sum\limits_{{ik} \in B} {\textit{NPV}{_{ik}} \cdot \textit{PoS}_{ik}^g \cdot {x_{ik}}}  + \sum\limits_{{jk} \in C} {\textit{NPV}{_{jk}} \cdot \textit{PoS}_{jk}^e \cdot {x_{jk}}} }}{{\sum\limits_{{ik} \in B} {{x_{ik}}}  + \sum\limits_{{jk} \in C} {{x_{jk}}} }}.$

\subsubsection{Constraints}

\textit{a}. The number of selected preliminary prospected wells and appraisal wells must equal the required number of such wells for the exploration deployment:
\begin{flalign}\label{18}
	&\sum\limits_{ik \in B} {\textit{Drill}_{ik}^{\textit{tra}} \cdot {x_{ik}}}  + \sum\limits_{jk \in C} {\textit{Drill}_{jk}^{\textit{app}} \cdot {x_{jk}}}  - \textit{Tot}{_{\textit{wells}}} = 0,&
\end{flalign}
where $ \textit{Drill}_{ik}^{\textit{tra}} $ denotes the number of preliminary prospected wells in the trap project, and $ \textit{Drill}_{jk}^{\textit{app}} $ denotes the number of appraisal wells in the appraisal projects.

\textit{b}. Pre-OR shall not be less than the lower bound of prospective oil reserves specified in the exploration plan:
\begin{flalign}\label{19}
	&\sum\limits_{{ik} \in B} \textit{Pred}_{ik}^{oi} \cdot {x_{ik}} - \textit{Pred}_{lb}^{oi} \ge 0,&
\end{flalign}

\textit{c}. Pre-GR shall not be less than the lower bound of prospective gas reserves specified in the exploration plan:
\begin{flalign}\label{20}
	&\sum\limits_{{ik} \in B} \textit{Pred}_{ik}^{ga} \cdot {x_{ik}} - \textit{Pred}_{lb}^{ga} \ge 0 ,&
\end{flalign}

\textit{d}. COR shall not be less than the lower bound of contingent oil reserves specified in the exploration plan:
\begin{flalign}\label{21}
	& \sum\limits_{{jk} \in C} \textit{Cont}_{jk}^{oi} \cdot {x_{jk}} - \textit{Cont}_{lb}^{oi} \ge 0,&
\end{flalign}

\textit{e}. CGR shall not be less than the lower bound of contingent gas reserves specified in the exploration plan:
\begin{flalign}\label{22}
	&\sum\limits_{{jk} \in C} \textit{Cont}_{jk}^{ga} \cdot {x_{jk}} - \textit{Cont}_{lb}^{ga} \ge 0 ,&
\end{flalign}

\textit{f}. Pro-OR shall not be less than the lower bound of proved oil reserves specified in the exploration plan:
\begin{flalign}\label{23}
	&\sum\limits_{{jk} \in C} \textit{Prov}_{jk}^{oi} \cdot {x_{jk}} - \textit{Prov}_{lb}^{oi} \ge 0,&
\end{flalign}

\textit{g}. Pro-GR shall not be less than the lower bound of proved gas reserves specified in the exploration plan:
\begin{flalign}\label{24}
	&\sum\limits_{{jk} \in C} \textit{Prov}_{jk}^{ga} \cdot {x_{jk}} - \textit{Prov}_{lb}^{ga} \ge 0 ,&
\end{flalign}

\textit{h}. The overall probability of success shall not be less than the average of the geological and economic probabilities of success:
\begin{flalign}\label{25}
	&\begin{array}{l}
		\frac{{\sum\limits_{ik \in B} {\textit{PoS}_{ik}^g\cdot \textit{Drill}_{ik}^{\textit{tra}}\cdot{x_{ik}}}  + \sum\limits_{jk \in C} {\textit{PoS}_{jk}^e\cdot \textit{Drill}_{jk}^{\textit{app}}\cdot{x_{jk}}} }}{{\sum\limits_{ik \in B} {\textit{Drill}_{ik}^{\textit{tra}}\cdot{x_{ik}}}  + \sum\limits_{jk \in C} {\textit{Drill}_{jk}^{\textit{app}}\cdot{x_{jk}}} }} - \textit{Drill}_{lb} \ge 0,
	\end{array}&
\end{flalign}

\textit{i}. Within the selected drilling portfolio, the number of projects with relatively low geological or economic probability of success (hereafter “success”) shall not exceed a predefined threshold:
\begin{flalign}\label{26}
	&\left\{ {\alpha :\;\textit{success}{_\alpha } < \textit{Thre}_{well},\;{z_\alpha } = 1} \right\} - {L_{\textit{ub} }} \le 0,&
\end{flalign}

\textit{j}. The total drilling cost for trap projects shall not exceed the company’s budget allocated for such projects:
\begin{flalign}\label{27}
	&\sum\limits_{{ik} \in B} {\textit{Cost}_{ik}^{\textit{tra}} \cdot {x_{ik}}}  - \textit{Cost}_{\textit{ub}}^{\textit{tra}} \le 0,&
\end{flalign}

\textit{k}. The total drilling cost for appraisal projects shall not exceed the company’s budget allocated for such projects:
\begin{flalign}\label{28}
	&\sum\limits_{{jk} \in C} {\textit{Cost}_{jk}^{\textit{app}} \cdot {x_{jk}}}  - \textit{Cost}_{\textit{ub}}^{\textit{app}} \le 0, &
\end{flalign}

\textit{l}. The optimized drilling portfolio should ensure that each region maintains a certain allocation:
\begin{flalign}\label{29}
	&N_{\textit{tra}}^{\left( r \right)} - \sum\limits_{{ik} \in {B_r}} {{x_{ik}}} \le 0,&
\end{flalign}
\begin{flalign}\label{30}
	&N_{\textit{app}}^{\left( r \right)} - \sum\limits_{{jk} \in {C_r}} {{x_{jk}}} \le 0,&
\end{flalign}

Therefore, the above formulation recasts the drilling deployment problem as a portfolio selection problem and, grounded in portfolio theory, models a multi-objective optimization model that jointly characterizes expected return and portfolio-level volatility under constraints on GPoS, three types of O\&G reserves, well count, investment budget, EPoS, and region. This multi-dimensional formulation enhances the robustness of the model.

\section{Operator-enhanced NSGA-II (OE-NSGA-II)}
\begin{figure*}[htbp]
	\centering
	\includegraphics[width=\linewidth]{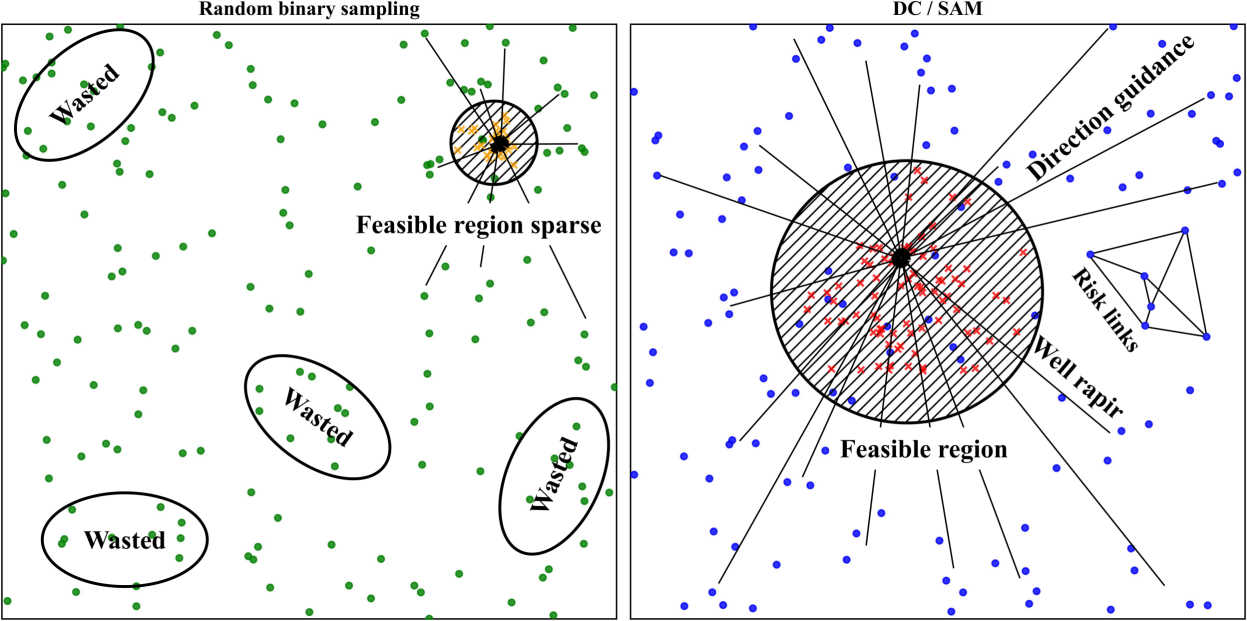} 
	\caption{Effect of initialization on feasible density and subsequent operator learning.}
	\label{fig3}
\end{figure*}

The NSGA-II algorithm combines fast non-dominated sorting with crowding-distance preservation and has demonstrated strong robustness in industrial multi-objective optimization tasks. Random-sampling–based initialization and operators are ill-suited to the highly constrained multi-objective model proposed herein, for two principal reasons: (1) The model enforces strong constraints—such as mandatory projects and fixed well counts—which render the feasible region sparse. Random initialization markedly reduces the proportion of feasible individuals, incurring substantial wasted evaluations and degrading computational efficiency. (2) The objective is defined at the portfolio level, capturing the return–risk trade-off and exhibiting strong nonseparability; random individuals fail to encode the risk-correlation structure of the portfolio. Consequently, during crossover and mutation, the information gain and useful learning signals transmitted from parents to offspring are limited, impairing convergence quality and the resulting Pareto front, shown in Fig. \ref{fig3}.

To address this limitation, we enhance the internal operators of the algorithm to better accommodate equality-constrained structures-the directional cross operator (DC), guided by parental decision bits, and a structure-aware mutation operator (SAM) with bit-level collisions are introduced: this strategy avoids modifying the objective functions while preserving non-dominated sorting and population diversity, thereby improving solution quality.

\subsection{The directional cross operator (DC) guided by parental decision bits}
To mitigate the adverse effects of random sampling on crossover and mutation, a directional crossover operator guided by parental decision bits is developed. The operator exploits directional information encoded in the parents’ key decision bits to steer recombination, thereby improving information transfer and strengthening directional search in the offspring.
\mysubtitle{1}{Normalization} 
In our model, the decision vector is $ z \in {\left\{ {0,1} \right\}^n} $, and each candidate project has expected return $ {g_i} > 0 $. Constraints include well-count requirements $ {w_i} \in\mathbb{Z}_{\ge0}^n $ (with $ {w_i} = 0 $ allowed for "reserve-providing appraisal wells"), a mandatory mask $ {m_i} \in \left\{ {0,1} \right\} $, and a region label $ r\left( i \right) \in R $. Define the selected set $ S\left( x \right) = \left\{ {i|{x_i} = 1} \right\} $,size $ n\left( x \right) = \left| {S\left( x \right)} \right| $, and mean $ \mu \left( x \right) = \frac{1}{{n\left( x \right)}}\sum\nolimits_{i \in S\left( x \right)} {{g_i}} $. We maximize total return $ \sum\nolimits_{i \in S\left( x \right)} {{g_i}} $ and minimize risk defined as the standard deviation of selected returns:
\begin{flalign}\label{31}
	& \begin{array}{l}
		R\left( x \right) = \sqrt {M\left( x \right)} ,\\
		M\left( x \right) = \sum\nolimits_{i \in S\left( x \right)} {{{\left( {{g_i} - \mu \left( x \right)} \right)}^2}}. 
	\end{array}&
\end{flalign}

To update portfolio-level risk, the second moment $ M $ is maintained using Welford’s method \cite{bib44}. Let the current statistics be $ \left( {n,\;\mu ,\;{M}} \right) $, denoting the sample size, mean, and second moment, respectively. 

Incremental update when adding a project with value $ g $: Let $ \delta  = g - \mu ,\;{n^ + } = n + 1,\;{\mu^ + } = {\delta  \mathord{\left/
		{\vphantom {\delta  {n^ + }}} \right.
		\kern-\nulldelimiterspace} {n^ + }} $ and $ \;{\delta^ +}  = g - {\mu^ + }$. Then $ \;{n^ + } = n + 1,\quad {\mu^ + } = \mu + \frac{g-\mu}{n^ + },\quad {M^ +} = M + (g-\mu)(g- {\mu^ + })\; $. 

Incremental update when removing a project with value $ g $: Let $ \delta  = g - \mu ,\;{n^ - } = n - 1,\;{\mu^ - } = {\delta  \mathord{\left/
		{\vphantom {\delta  {n'}}} \right.
		\kern-\nulldelimiterspace} {n^ -}} $. Then $ \;{n^ - } = n-1,\quad {\mu^ - } = \mu - \frac{g-\mu}{n^ -},\quad{ M^ -} = M - (g-\mu)(g-{\mu^ - })\; $. Risk increment definition: If risk is measured by the sample variance $ \Delta M_i^{0 \to 1} = {M^ - } - M $, $ \Delta M_i^{1 \to 0} = {M^ - } - M $. where $ \Delta M_i^{0 \to 1} $ and $ \Delta M_i^{1 \to 0} $ denote the increment in the risk measure when flipping the $ i^{th} $ bit from $ 0 \to 1 $ or $ 1 \to 0 $, respectively.

To balance expected return against risk, a directional parameter $ \rho  \sim Beta\left( {\alpha ,\;\alpha } \right),\;\alpha  > 0 $ is introduced. As $ \rho  \to 1 $, the preference shifts toward return; as $ \rho  \to 0 $, it emphasizes risk reduction. For any vector $ u $ on the candidate index set $ J $, min–max normalization is applied:
\begin{flalign}\label{32}
	&{\hat u_i} = \frac{{{u_i} - {{\min }_{j \in J}}{u_j}}}{{{{\max }_{j \in J}}{u_j} - {{\min }_{j \in J}}{u_j} + \varepsilon }}. &
\end{flalign}

To enforce the quota $ \sum\nolimits_{i:r\left( i \right) = c} {{x_i} \ge {q_c}} $, a shortfall-guided bias is introduced for each candidate $ i $: 
\begin{flalign}\label{33}
	&\begin{array}{l}
		{b_i}\left( x \right) = k.\max \left\{ {0,\;{q_{r\left( i \right)}} - {C_{r\left( i \right)}}\left( x \right)} \right\}, \;k \ge 0 ,\\
	{C_c}\left( x \right) = \sum\nolimits_{j:r\left( j \right) = c} {{x_j}},
     \end{array}&
\end{flalign}
where $ r\left( i \right) $ maps project $ i $ to region $ c $, $ {q_c} \ge 0 $ is the minimum quota for region $ c $, $ {C_c}\left( x \right) $ counts selected projects in region $ c $, and $ k $ tunes the guidance strength. When a region is underfilled $ {C_c}\left( x \right) < {q_c} $, projects in that region receive a positive $ {b_i}\left( x \right) $, thereby biasing selection toward 1 (via probabilities, scores, or operator priors) and accelerating convergence to the feasible set.

For any current solution $ x $, directional preference scores are defined for the $ i^{th} $ bit taking values 0 or 1: 
\begin{flalign}\label{34}
	&\begin{array}{l}
		{\rm{Directio}}{{\rm{n}}_i}\left( 1 \right) = \rho {{\hat g}_i} - (1 - \rho )\gamma \Delta {M^{0 \to 1}} + {b_i}\left( x \right), \\
	{\rm{Directio}}{{\rm{n}}_i}\left( 0 \right) =  - \rho {{\hat g}_i} - (1 - \rho )\gamma \Delta {M^{1 \to 0}},
    \end{array}&
\end{flalign}
where $ {{\hat g}_i} $ is the min–max normalized expected return over the candidate index set $ \gamma  > 0 $ is the risk weight, and $ {b_i}\left( x \right) $ is the regional shortfall-guidance term. Decision rule: If $ {\rm{Directio}}{{\rm{n}}_i}\left( 1 \right) \ge {\rm{Directio}}{{\rm{n}}_i}\left( 0 \right) $, set $ {x_i} = 1 $; otherwise set $ {x_i} = 0 $.
\mysubtitle{2}{Scope and initialization}
Given two binary parents $ A,\;B \in {\left\{ {0,1} \right\}^n} $, DC recombines only the differing loci $ D = \left\{ {i|\;{A_i} \ne {B_i}} \right\} $ to avoid unnecessary perturbations and preserve parental consensus. \textit{offspring}${_1} $ is initialized from $ A $; its selection statistics $ \left( {n,\;\mu ,\;{M}} \right) $ are computed once using $ R\left( x \right) $ and $ M\left( x \right) $. The complementary \textit{offspring}${_2} $ mirrors the procedure using $ B $.
\mysubtitle{3}{Directional decision on differing loci}
For each $ i \in D $, let the min–max normalized return and incremental risk be $ {\hat g_i} $ and $ \Delta {M_i} $; let $ \rho  \sim Beta\left( {\alpha ,\;\alpha } \right) $ be the preference parameter and $ {b_i}\left( x \right) $ the optional regional soft bias. The per-locus directional scores are $ {\rm{Directio}}{{\rm{n}}_i}\left( 1 \right) $ and $ {\rm{Directio}}{{\rm{n}}_i}\left( 0 \right) $. The locus decision is \textit{offspring}$ _{1,i} = \mathbb{I} \left\{ {{\rm{Directio}}{{\rm{n}}_i}\left( 1 \right) \ge {\rm{Directio}}{{\rm{n}}_i}\left( 0 \right)} \right\} $. Finally, must-select constraints are enforced: if $ {m_i} = 1 $ then \textit{offspring}${_{1,i}} \leftarrow 1 $.
\mysubtitle{4}{Greedy well-count repair}
Let the well target be $ {W^*} $ and define the well deficit as $ \Delta W = W^{*} - \sum_{i} w_i\, \textit{offspring}_{1,i} $, DC repairs only actionable loci $ C = \{ i \in D|{w_i} \ge 1\} $ using a unit-well benefit ranking.

(Case) $ \Delta W > 0 $ (need to add wells): consider $ {C_{add}} = \left\{ {i \in C|\textit{offspring}{_{1,i}} = 0} \right\} $ and rank by $ \frac{{\rho {{\hat g}_i} - \left( {1 - \rho } \right)\gamma \Delta M_i^{0 \to 1} + {b_i}\left( x \right)}}{{\max \left\{ {1,{w_i}} \right\}}} $, switching bits to 1 until  $ {W^*} $ is met. 

\begin{figure*}[htbp]
	\centering
	\includegraphics[width=0.8\linewidth]{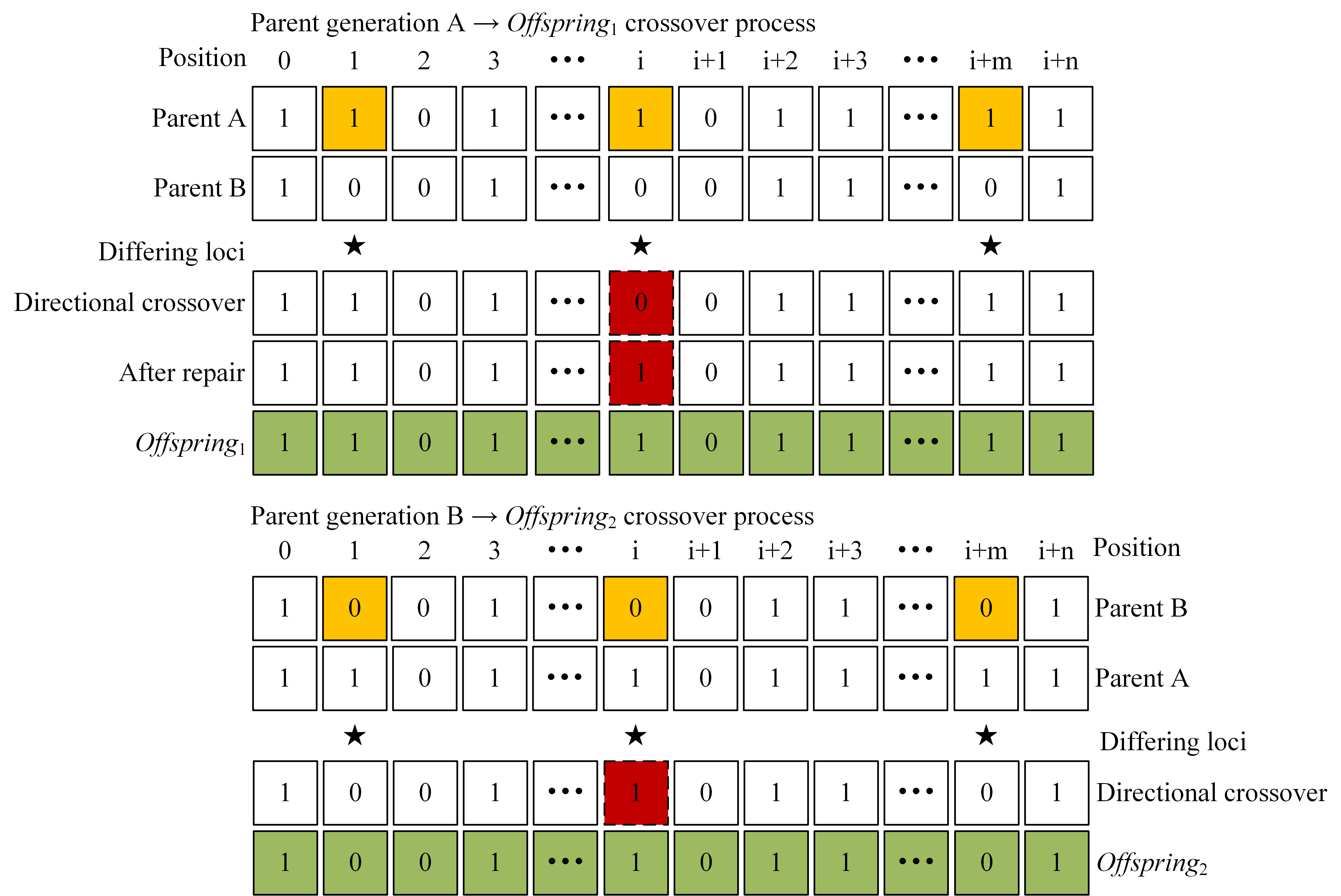} 
	\caption{Schematic diagram of the directional cross operator guided by parental decision bits.}
	\label{fig4}
\end{figure*}

(Case) $ \Delta W < 0 $ (need to add wells): consider $ {C_{rem}} = \left\{ {i \in C|\textit{offspring}{_{1,i}} = 1} \right\} $ and rank by $ \frac{{\rho {{\hat g}_i} - \left( {1 - \rho } \right)\gamma \Delta M_i^{1 \to 0} + {b_i}\left( x \right)}}{{\max \left\{ {1,{w_i}} \right\}}} $, switching bits to 0 until  $ {W^*} $ is met.  A schematic of the process is shown in Fig. \ref{fig4}. The detailed procedure of the crossover operator is provided in algorithmic \ref{alg:1}.

\begin{breakablealgorithm}
	\caption{DC: The directional cross operator}
	\label{alg:1}
	\begin{algorithmic}[1]
		\Require Parents $A,B\in\{0,1\}^n$; $g\in\mathbb{R}^n$, $w\in\mathbb{Z}_{\ge0}^n$, $m\in\{0,1\}^n$; $r(\cdot)$, $\{q_c\}$; $\alpha>0$, $\gamma>0$, $k\ge 0$; $W^\star$
		\Ensure $\textit{offspring}_1,\textit{offspring}_2\in\{0,1\}^n$
		
		\State $D \leftarrow \{i~|~A_i\neq B_i\}$; \quad $\rho \sim \mathrm{Beta}(\alpha,\alpha)$ 
		\State $\textit{offspring}_1 \leftarrow A$;\quad $(n,\mu,M)\leftarrow (\textit{offspring}_1,g)$
		\If{$|D|>0$}
		\State $\hat g_D \leftarrow \text{MINMAX}\big(g[D]\big)$
		\State $d^{+}[i\in D] \leftarrow \Delta {M^{0 \to 1}}(n,\mu,M,g_i)$;\quad $d^{-}[i\in D] \leftarrow \Delta {M^{1 \to 0}}(n,\mu,M,g_i)$
		\State $\hat d^{+}\leftarrow \text{MINMAX}(d^{+})$;\quad $\hat d^{-}\leftarrow \text{MINMAX}(d^{-})$
		\For{$i\in D$} 
		\State $b_i \leftarrow (i,\textit{offspring}_1,r(\cdot),q,k)$
		\State $\text{Direction(1)} \leftarrow \rho\,\hat g_D[i] - (1-\rho)\,\gamma\,\hat d^{+}[i] + b_i$
		\State $\text{Direction(0)} \leftarrow -\rho\,\hat g_D[i] - (1-\rho)\,\gamma\,\hat d^{-}[i]$
		\State $\textit{offspring}_1[i] \leftarrow \mathbb{I}\{\text{Direction(1)}\ge \text{Direction(0)}\}$
		\EndFor
		\EndIf
		
		\For{$i=1$ to $n$}
		\If{$m_i=1$}
		\State $\textit{offspring}_1[i]\leftarrow 1$ 
		\EndIf
		\EndFor
		
		\If{$W^\star$ is not None} 
		\State $\Delta W \leftarrow W^\star - \sum_i w_i\,\textit{offspring}_1[i]$
		\If{$\Delta W \neq 0$}
		\State $C \leftarrow \{i\in D~|~w_i\ge 1\}$
		\If{$\Delta W>0$}
		\State $\text{benefit}[i\in C] \leftarrow \rho\,\hat g(i) - (1-\rho)\,\gamma\,\big(d^{+}(i)\big) + (i,\textit{offspring}_1,\cdot)$
		\State $\text{order} \leftarrow \text{argsort}_{\downarrow}\big(\text{benefit}[i]/\max\{1,w_i\}\big)$
		\For{$i$ in \text{order}}
		\If{$\Delta W\le 0$} 
		\State \textbf{break}
		\EndIf
		\If{$\textit{offspring}_1[i]=0$}
		\State $\textit{offspring}_1[i]\leftarrow 1$
		\EndIf
		\EndFor
		\Else
		\State $\text{benefit}[i\in C] \leftarrow \rho\,\hat g(i) - (1-\rho)\,\gamma\,\big(d^{-}(i)\big) + (i,\textit{offspring}_1,\cdot)$
		\State $\text{order} \leftarrow \text{argsort}_{\uparrow}\big(\text{benefit}[i]/\max\{1,w_i\}\big)$
		\For{$i$ in \text{order}}
		\If{$\Delta W\ge 0$} 
		\State \textbf{break}
		\EndIf
		\If{$\textit{offspring}_1[i]=1$}
		\State $\textit{offspring}_1[i]\leftarrow 0$
		\EndIf
		\EndFor
		\EndIf
		\EndIf
		\EndIf
		\State Repeat Lines 1--42 with $A\leftarrow B$ and $\rho\leftarrow 1-\rho$ to obtain $\textit{offspring}_2$
		
		\State \Return $\textit{offspring}_1,\textit{offspring}_2$
	\end{algorithmic}
\end{breakablealgorithm}

\subsection{Structure-aware mutation (SAM) operator with bit-level collisions}
\begin{figure*}[htbp]
	\centering
	\includegraphics[width=0.8\linewidth]{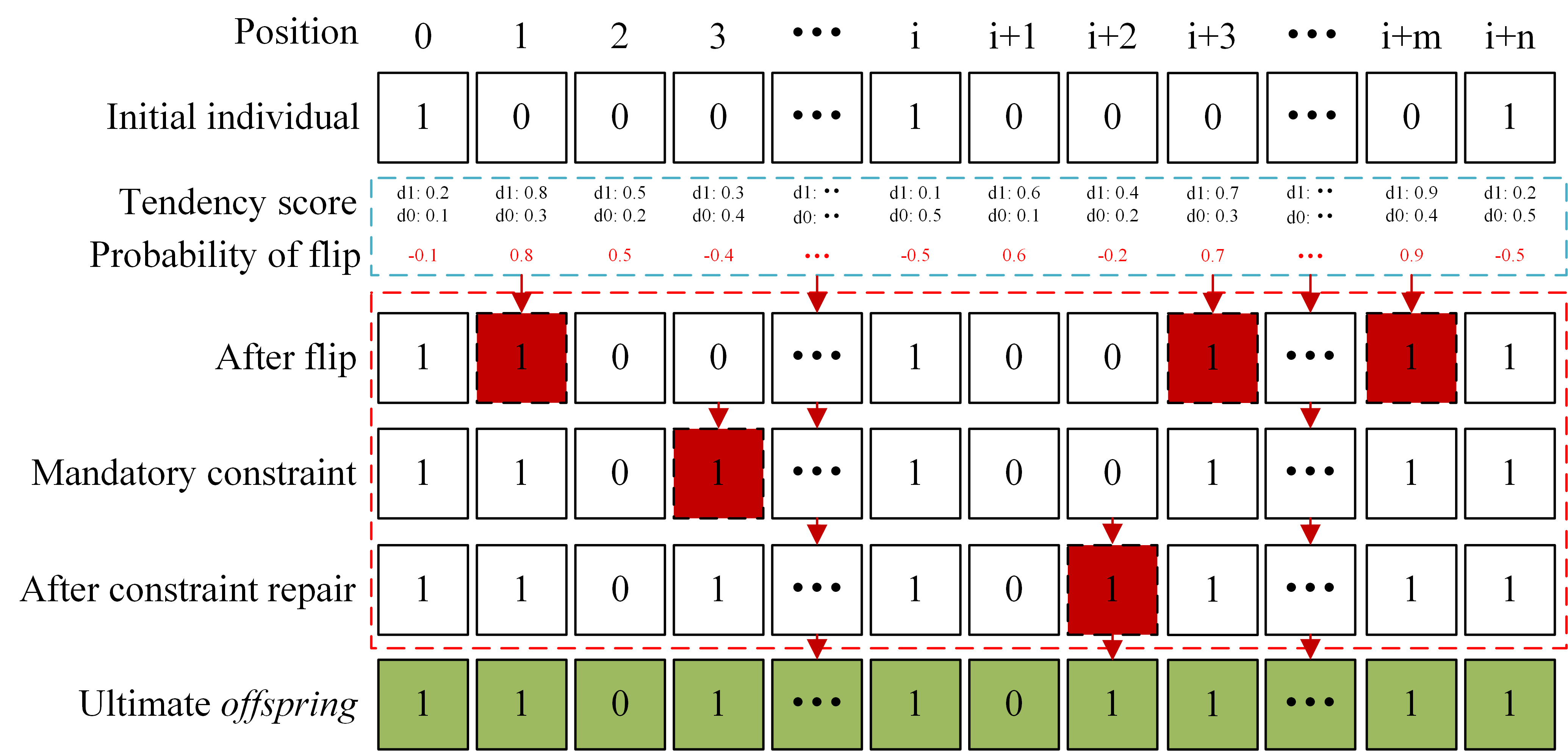} 
	\caption{Schematic diagram of a structure-aware mutation operator with bit-level collisions.}
	\label{fig5}
\end{figure*}

The improved crossover operator can partially guide information transfer across decision bits, but the standard bit-flip mutation still flips bits “blindly” with independent probabilities, without sensing the nonseparability of the objective function and thus misaligned with the model structure. To address this, we propose a structure-aware mutation (SMA) operator with bit-level collisions that uses probabilistic sampling combined with feasibility repair under constraints to perform directed search within the feasible region toward a better return–risk balance, thereby strengthening the exploration of high-quality Pareto solutions.
\mysubtitle{1}{Budgeted directional mutation}
Given an individual $ x \in {\left\{ {0,\;1} \right\}^n} $, sample $ \rho  \sim Beta\left( {\alpha ,\;\alpha } \right) $ set a mutation budget $ L = \max \left\{ {{L_{\min }},\;\left\lceil \beta  \right\rceil } \right\},\;0 < \beta  \ll 1 $, to cap the number of flips per individual. Compute $ \left( {n,\;\mu ,\;{M}} \right) $ once for $ x $.
\mysubtitle{2}{Scoring and selecting flips}
Compute in $ O\left( n \right) $ the incremental risk vectors $ \Delta {M^{1 \to 0}} $ and  $ \Delta {M^{0 \to 1}} $, normalize them and form the directional scores. Define a unified flip gain:
\begin{flalign}\label{35}
	& \Phi  = \left\{ \begin{array}{l}
		{\rm{Directio}}{{\rm{n}}_i}\left( 1 \right),\;{x_i} = 0 \\ 
		{\rm{-Directio}}{{\rm{n}}_i}\left( 0 \right),\;{x_i} = 1 \\ 
	\end{array} \right. &
\end{flalign}
which measures the desirability of changing the current bit. Select the top-L indices by $ \Phi $ (descending) and flip $ {x_i} \leftarrow 1 - {x_i} $.
\mysubtitle{3}{Constraint enforcement and well-count repair}
Enforce must-select: if $ {m_i} = 1 $, set $ {x_i} \leftarrow 1 $. If a hard well-count target $ {W^*} $ is active, compute $ \Delta W = {W^*} - \sum\nolimits_i {{w_i}{x_i}} $ and perform the same greedy unit-well repair as in DC, now over the natural add/remove candidate sets. A schematic of the process is shown in Fig. \ref{fig5}. The detailed procedure of the mutation operator is provided in algorithmic \ref{alg:2}.

\begin{breakablealgorithm}
	\caption{SAM: Structure-aware mutation operator}
	\label{alg:2}
	
	\begin{algorithmic}[1]
		\Require $x\in\{0,1\}^n$; $g,w,m,r(\cdot),q$; $\alpha>0$, $\gamma>0$, $k\ge 0$; budget ratio $\beta\in(0,1)$, lower bound $L_{\min}$; $W^\star$ 
		\Ensure Mutated $x'\in\{0,1\}^n$
		
		\State $x' \leftarrow x$;\quad $\rho \sim \mathrm{Beta}(\alpha,\alpha)$;\quad $(n,\mu,M)\leftarrow (x',g)$
		\State $L \leftarrow \max\big\{L_{\min},\lceil \beta n\rceil\big\}$
		
		\State $\hat g \leftarrow \text{MINMAX}(g)$
		\State $d^{+}[i] \leftarrow \Delta {M^{0 \to 1}}(n,\mu,M,g_i)$;
		\State $d^{-}[i] \leftarrow \Delta {M^{1 \to 0}}(n,\mu,M,g_i)$
		\State $\hat d^{+}\leftarrow \text{MINMAX}(d^{+})$;\quad $\hat d^{-}\leftarrow \text{MINMAX}(d^{-})$
		\State $b[i] \leftarrow (i,x',r(\cdot),q,k)$ for all $i$
		
		\State $\text{Direction(1)}[i] \leftarrow \rho\,\hat g[i] - (1-\rho)\,\gamma\,\hat d^{+}[i] + b[i]$
		\State $\text{Direction(0)}[i] \leftarrow -\rho\,\hat g[i] - (1-\rho)\,\gamma\,\hat d^{-}[i]$
		
		\State $C_0 \leftarrow \{i~|~x'_i=0\}$ with $\text{Score}[i]=\text{Direction(1)}[i]$
		\State $C_1 \leftarrow \{i~|~x'_i=1\}$ with $\text{Score}[i]=-\text{Direction(0)}[i]$
		\State $C \leftarrow C_0 \cup C_1$;
		\For{$i\in T$}
		\State $x'_i \leftarrow 1-x'_i$
		\EndFor
		
		\For{$i=1$ to $n$}
		\If{$m_i=1$}
		\State $x'_i\leftarrow 1$ 
		\EndIf
		\EndFor
		
		\If{$W^\star$ is not None} 
		\State $\Delta W \leftarrow W^\star - \sum_i w_i\,x'_i$
		\If{$\Delta W>0$}
		\State $A \leftarrow \{i~|~x'_i=0,\,w_i\ge 1\}$
		\State $\text{benefit}[i\in A] \leftarrow \rho\,\hat g[i] - (1-\rho)\,\gamma\,\text{MINMAX}\big(d^{+}[i]\big)$
		\State $\text{order} \leftarrow \text{argsort}_{\downarrow}\big(\text{benefit}[i]/\max\{1,w_i\}\big)$
		\For{$i$ in \text{order}}
		\If{$\Delta W\le 0$}
		\State \textbf{break}
		\EndIf
		\State $x'_i \leftarrow 1$;\quad $\Delta W \leftarrow \Delta W - w_i$
		\EndFor
		\EndIf
		\If{$\Delta W<0$}
		\State $R \leftarrow \{i~|~x'_i=1,\,w_i\ge 1\}$
		\State $\text{benefit}[i\in R] \leftarrow \rho\,\hat g[i] - (1-\rho)\,\gamma\,\text{MINMAX}\big(d^{-}[i]\big)$
		\State $\text{order} \leftarrow \text{argsort}_{\uparrow}\big(\text{benefit}[i]/\max\{1,w_i\}\big)$
		\For{$i$ in \text{order}}
		\If{$\Delta W\ge 0$}
		\State \textbf{break}
		\EndIf
		\State $x'_i \leftarrow 0$;\quad $\Delta W \leftarrow \Delta W + w_i$
		\EndFor
		\EndIf
		\EndIf
		
		\State \Return $x'$
	\end{algorithmic}
\end{breakablealgorithm}

\section{Experimental results}
This section comprises the following components: data sources, multi-objective algorithm evaluation metrics, algorithm parameter settings, method validation, method application, and the resulting exploration investment schemes. All numerical experiments were conducted on a workstation equipped with a 12th Gen Intel Core i5-12490F (3.00 GHz) processor and 16 GB RAM.

\subsection{Data}
\begin{table*}[htbp]
	\caption{Prospect investment lists of trap projects and core parameters.}\label{Tab1}
	\centering
	\resizebox{\textwidth}{!}{%
		\begin{tabular}{lcccccccc}
			\hline
			Region & Trap projects & Pre-OR ($10^4$ t) & Pre-GR ($10^8$ $m^3$) & Cost ($10^4$ CNY) & NPV ($10^4$ CNY) & GPoS (\%) & Well count & mandatory \\
			\midrule
			E & QL3  & 38.80  & 3.70   & 3087.00   & 13515.00   & 0.53  & 1.00  & 0.00   \\ 
			E & YQZY1  & 138.20  & 0.00   & 3112.00   & 7335.00   & 0.26  & 1.00  & 0.00   \\ 
			B & YQX12  & 285.80  & 0.00   & 4017.00   & 21832.00   & 0.76  & 1.00  & 0.00   \\ 
			B & YQX14  & 2443.79  & 0.00   & 4024.00   & 10851.00   & 0.76  & 1.00  & 0.00   \\ 
			D & BST1  & 19.40  & 1.60   & 7355.00   & 406.00   & 0.16  & 1.00  & 1.00   \\ 
			C & Z1A1  & 30.00  & 10.00   & 8593.00   & 5862.00   & 0.69  & 1.00  & 0.00   \\ 
			E & XH7  & 110.40  & 125.20   & 8663.00   & 9652.00   & 0.82  & 1.00  & 0.00   \\ 
			A & SB42  & 610.00  & 0.00   & 10048.00   & 57159.00   & 0.45  & 1.00  & 0.00   \\ 
			A  & SB6 & 476.20  & 129.20  & 10047.00   & 69826.00  & 0.93  & 1.00  & 0.00   \\ 
			A  & SB14 & 0.00  & 298.60  & 10096.00   & 26182.10  & 0.68  & 1.00  & 0.00   \\ 
			A  & SBL2 & 373.00  & 0.00  & 11548.00   & 26078.00  & 0.50  & 1.00  & 0.00   \\ 
			A  & SBL3 & 1669.14  & 0.00  & 11659.00   & 4855.00  & 0.82  & 1.00  & 0.00   \\ 
			A  & SB13F2 & 1188.60  & 0.00  & 11847.00   & 39206.06  & 0.32  & 1.00  & 0.00   \\ 
			A  & SB8N1 & 482.00  & 216.00  & 11864.00   & 36108.62  & 0.75  & 1.00  & 0.00   \\ 
			A  & SN3 & 553.20  & 371.00  & 11883.00   & 40047.17  & 0.47  & 1.00  & 0.00   \\ 
			A  & SB4 & 385.20  & 148.37  & 11866.00   & 41200.80  & 0.41  & 1.00  & 0.00   \\ 
			A  & SB10F2 & 235.80  & 101.72  & 10986.00   & 87075.12  & 0.66  & 1.00  & 0.00   \\ 
			A  & SB8F6 & 243.40  & 52.20  & 10807.00   & 15253.25  & 0.95  & 1.00  & 0.00   \\ 
			A  & SB9 & 877.80  & 0.00  & 11988.00   & 34282.65  & 0.38  & 1.00  & 0.00   \\ 
			A  & SB11F2 & 368.40  & 0.00  & 12053.00   & 30677.00  & 0.53  & 1.00  & 0.00   \\ 
			A  & SB10 & 285.60  & 116.40  & 12420.00   & 29443.82  & 0.51  & 1.00  & 0.00   \\ 
			A  & SB11F3 & 550.20  & 0.00  & 12234.00   & 32352.00  & 0.52  & 1.00  & 0.00   \\ 
			A  & SB8N2 & 415.20  & 189.00  & 12980.00   & 70645.00  & 0.67  & 1.00  & 0.00   \\ 
			D  & TS2 & 759.60  & 37.80  & 12985.00   & 73517.00  & 0.04  & 1.00  & 0.00   \\ 
			A  & SB8NY1 & 92.40  & 71.40  & 13587.00   & 25278.30  & 0.92  & 1.00  & 0.00   \\ 
			C & KL3 & 772.00  & 143.60  & 14168.00   & 60539.00  & 12461.00  & 1.00  & 0.00   \\ 
			D & TSX1 & 3247.40  & 159.80  & 14268.00   & 81576.00  & 12922.43  & 1.00  & 0.00   \\ 
			D & TSW1 & 1434.00  & 0.00  & 14829.00   & 22986.00  & 4433.14  & 1.00  & 0.00   \\ 
			A & SB2F1 & 1053.29  & 243.93  & 17663.00   & 29519.00  & 5200.34  & 1.00  & 0.00   \\ 
			C & K4X1 & 0.00  & 435.40  & 22563.00   & 18258.00  & 2240.33  & 1.00  & 1.00   \\ 
			\bottomrule 
	\end{tabular}}
\end{table*}
\begin{table*}[htbp]
	\caption{Prospect investment lists of appraisal projects and core parameters.}\label{Tab2}
	\centering
	\resizebox{\textwidth}{!}{%
		\begin{tabular}{lcccccccccc}
			\hline
			Region & Appraisal projects & COR ($10^4$ t) & CGR ($10^8$ $m^3$) & Pro-OR ($10^4$ t) & Pro-GR ($10^8$ $m^3$) & Cost ($10^4$ CNY) & NPV ($10^4$ CNY) & EPoS (\%) & Well count & Mandatory \\
			\midrule
			A & SB12X  & 0 & 460 & 0 & 0 & 0 & 83450  & 0.6 & 0 & 0  \\ 
			A & SB1X  & 0 & 0 & 180 & 100 & 8953 & 18344  & 0.609 & 2 & 0  \\ 
			B & TH10  & 169 & 0 & 0 & 0 & 0 & 11012  & 0.9025 & 0 & 0  \\ 
			B & TH125  & 0 & 0 & 986 & 0 & 0 & 19465  & 0.8 & 0 & 0  \\ 
			C & YB1  & 131.47 & 0 & 0 & 0 & 6890 & 4881  & 0.88 & 1 & 0  \\ 
			C & Z1  & 0 & 0 & 213 & 85 & 12060 & 17868  & 0.65 & 1 & 0  \\ 
			A & SZ41  & 1898.02 & 504.32 & 0 & 0 & 8500 & 99832  & 0.89 & 1 & 1  \\ 
			B & TH122  & 1982.64 & 0 & 0 & 0 & 0 & 1310  & 0.95 & 0 & 0  \\ 
			B & YQ516  & 1782.56 & 0 & 0 & 0 & 0 & 4800  & 0.86 & 0 & 0  \\ 
			D & AT28  & 106.43 & 0 & 0 & 0 & 0 & 1800  & 0.95 & 0 & 0  \\
			D & S81  & 54.24 & 0 & 0 & 0 & 0 & 1860  & 0.95 & 0 & 1  \\ 
			A & S9  & 108.46 & 0 & 0 & 0 & 0 & 65  & 0.95 & 0 & 2  \\ 
			A & SB10X  & 0 & 0 & 363.09 & 306.11 & 8700 & 1232  & 0.85 & 1 & 0  \\ 
			A & SZ412  & 0 & 0 & 980.74 & 193.64 & 9000 & 51000  & 0.89 & 1 & 0  \\ 
			B & YQX11  & 0 & 0 & 4034.12 & 0 & 4800 & 2365  & 0.82 & 1 & 0  \\ 
			\bottomrule 
	\end{tabular}}
\end{table*}
All data used in this study are derived from the company’s real exploration operations. To protect sensitive information, the dataset has been anonymized and sanitized in accordance with standard confidentiality practices. Since the prospect lists for 2023 and 2024 differ only in their candidate projects, we report a subset of the 2023 prospect list in Tables \ref{Tab1} and \ref{Tab2}. It is worth noting that key quantities such as GPoS, EPoS, the three types of O\&G reserves, and NPV in these tables are generated using the methodology outlined in Section 3.
\begin{table*}[htbp]
	\caption{Different algorithm parameter settings.}\label{Tab3}
	\centering
	\resizebox{\textwidth}{!}{%
		\begin{tabular}{lcc}
			\hline
			Algorithm & Parameters & value \\
			\midrule
			NSGA-II & Crossover probability, Mutation probability & $ prob = 0.9$, $ pm =0.05 $ \\
			NSGA-III & Crossover probability, Mutation probability, Number of reference points & $ prob = 0.9$, $ pm =0.05 $,  \textit{numRef} = \textit{numPop}  \\
			UNSGA-III & Crossover probability, Mutation probability, Reference directions & $ prob = 0.9$, $ pm =0.05 $, $ {n_ - }{\textit{points}} = 99  $ \\
			AGE-MOEA & Adaptive operator,  & $p \in \left[ {0.1,\;20} \right]$  \\
			AGE-MOEA2 & p-norm estimation,  & $p \in \left[ {0.1,\;20} \right]$  \\
			RVEA & Penalty function parameters, Reference vector adaptive frequency  & $ \alpha = 2 $, $ {\textit{adapt}_ - }{\textit{freq}} = 0.1  $  \\
			\textbf{OE-NSGA-II} & distribution shape parameter, soft-bias strength, risk weight, mutation budget & $\alpha  = 0.7$, $k  = 0.3$, $\gamma  = 1.3$, $\beta  = 0.05$ \\
			\bottomrule 
	\end{tabular}}
\end{table*}

\subsection{Performance indicators}
To assess the performance of OE-NSGA-II in this application, it is benchmarked against representative multi-objective evolutionary algorithms, including NSGA-II \cite{bib24}, NSGA-III \cite{bib45}, U-NSGA-III \cite{bib46}, AGE-MOEA \cite{bib47}, AGE-MOEA-II \cite{bib48}, and RVEA \cite{bib49}. Evaluation is based on three widely used metrics: (1) Hypervolume (HV) \cite{bib50}, (2) Inverted generational distance (IGD) \cite{bib51}, (3) Spacing \cite{bib52,bib53}, and (4) Set coverage (SC) \cite{bib54}.

HV quantifies the extent to which an approximation set covers the objective space. Given a reference point $ r $ (for minimization, chosen to be dominated by all solutions), HV is the dominated volume enclosed by the approximate Pareto set $ S $ and $ r $. If an approximation set $ A $ strictly dominates $ B $, then $ \mathrm{HV}\left( A \right) > \mathrm{HV}\left( B \right) $. Hence, HV captures both convergence (proximity to the true Pareto front) and diversity (spread along the front). The formal definition is as follows:
\begin{flalign}\label{36}
	& \mathrm{HV}\left( {S;\;r} \right) = \lambda \left( {{ \cup _{x \in S}}\left[ {{f_1}\left( x \right),\;{r_1}} \right] \times  \ldots  \times \left[ {{f_n}\left( x \right),\;{r_n}} \right]} \right). &
\end{flalign}

IGD captures both convergence and diversity. For each point $ y $ in the reference set $ \textit{PF}^* $ of the true Pareto front, its nearest neighbor $ x $ in the approximation set $ \textit{PF} $ is identified and their Euclidean distances in the objective space are averaged:
\begin{flalign}\label{37}
	& \mathrm{IGD} = \frac{1}{{\left| \textit{PF}^{*} \right|}}\sum\nolimits_{y \in \textit{PF}^{*}} {d\left( {y,\;x} \right)}, &
\end{flalign}
a smaller IGD indicates an approximation set that is closer to, and more evenly spread along, the true front. The reference set $ \textit{PF}^{*} $ is typically dense and uniform (often larger than $ \textit{PF} $ to enable comprehensive evaluation.
\begin{figure*}[htbp]
	\centering
	\includegraphics[width=\linewidth]{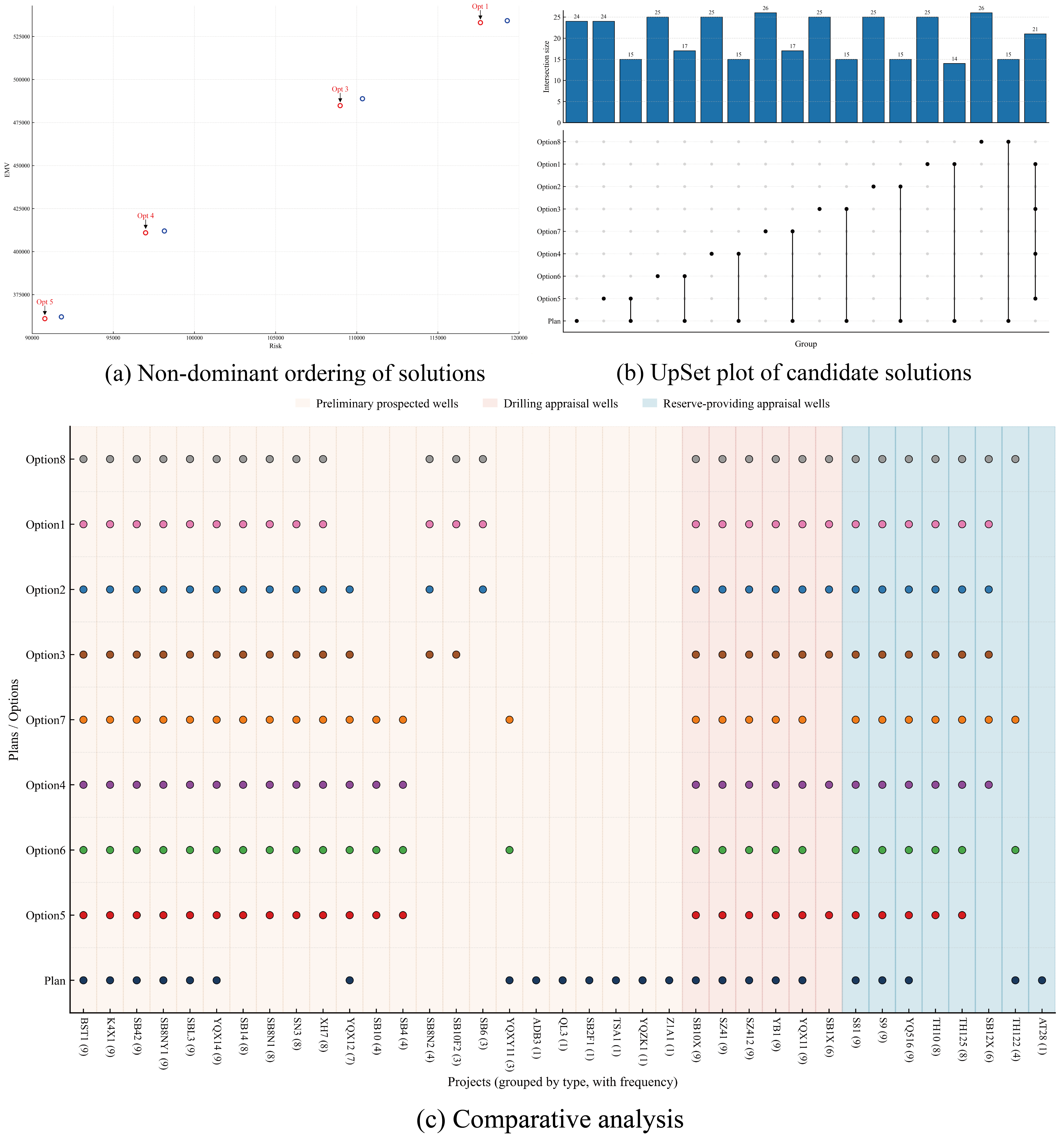} 
	\caption{Detailed comparison between selected projects and the exploration plan.}
	\label{fig6}
\end{figure*}
\begin{table*}[htbp]
	\caption{HV, IGD, and Spacing values obtained by all algorithms.}\label{Tab4}
	\centering
	\begin{tabular}{C{2.8cm}C{2cm}C{2cm}C{2cm}C{3cm}C{3cm}}
		\midrule
		Algorithm & HV  & IGD  & Spacing  & SC (OE-NSGA-II, X) & SC (X, OE-NSGA-II)\\ 
		\hline
		OE-NSGA-II & 5.4877E+09  & 2.6342E+02  & 2.3363E+02  & Nan  & Nan\\ 
		NSGA-III & 5.3562E+09  & 1.1253E+03  & 1.2531E+03  & 4.7059E-01  & 1.9355E-01 \\ 
		UNSGA-III & 5.3977E+09  & 6.3544E+02  & 4.4352E+03  & 6.5909E-01 & 6.4516E-02 \\ 
		NSGA-II & 4.0285E+09  & 1.2502E+03  & 1.3875E+03  & 3.4694E-01 &  3.2484E-01\\ 
		RVEA & 4.4320E+09  & 1.8523E+03  & 4.6250E+03  & 1.0000E+00 & 1.3000E-01\\ 
		AGE-MOEA-II & 4.3365E+09  & 1.4896E+03  & 5.9235E+02 & 4.8387E-01 & 2.6263E-01 \\ 
		AGE-MOEA & 4.3588E+09  & 1.4966E+03  & 3.8630E+03  & 3.900E-01 & 2.2581E-01 \\ 
		\bottomrule 
	\end{tabular}
\end{table*}

Spacing measures the uniformity of an approximate Pareto set in objective space: a smaller value indicates a more even spread. Let $ N $ be the number of solutions and $ {f^{\left( i \right)}} = \left( {f_1^{\left( i \right)},\;f_2^{\left( i \right)},\; \ldots ,\;f_M^{\left( i \right)}} \right) $ the objective vector of solution $ i $. Define the nearest-neighbor Manhattan distance $ {d_i} = \mathop {\min }\limits_{j \ne i} \sum\nolimits_{m = 1}^M {\left| {f_m^{\left( i \right)} - f_m^{\left( j \right)}} \right|} $, and $ \hat d = \frac{1}{N}\sum\nolimits_{i = 1}^N {{d_i}}  $. The spacing is:
\begin{flalign}\label{38}
	& \text{Spacing} = \sqrt {\frac{1}{{N - 1}}\sum\nolimits_{i = 1}^N {{{\left( {{d_i} - \hat d} \right)}^2}} }. &
\end{flalign}

SC evaluates the dominance relationship between two Pareto fronts. Let $ A $ and $ B $ denote two approximation sets, the coverage of $ B $ by $ A $ is defined as:
\begin{flalign}\label{39}
	& \mathrm{SC}\left( {A,\;B} \right) = \frac{{\left| {\left\{ {b \in B|\;\exists a \in A:\;a \prec b} \right\}} \right|}}{{\left| B \right|}}, &
\end{flalign}
where $ a \prec b $ denotes that $ a $ is no worse than $ b $ in all objectives and strictly better in at least one (for minimization). $ \mathrm{SC}\left( {A,\;B} \right) \in \left[ {0,\;1} \right] $ and is asymmetric; it is customary to report both $ \mathrm{SC}\left( {A,\;B} \right) $ and $ \mathrm{SC}\left( {B,\;A} \right) $. If $ \mathrm{SC}\left( {A,\;B} \right) = 1 $ and $ \mathrm{SC}\left( {B,\;a} \right) = 0 $, then $ A $ fully dominates $ B $, similar values indicate comparable or mutually non-dominating performance.

\subsection{Parameter settings for OE-NSGA-II}
In the practical tests, minor tuning of key OE-NSGA-II parameters was performed to ensure effectiveness: the maximum number of iterations was set to 500 and the population size to 100; the distribution shape parameter was $\alpha  = 0.7$, the soft-bias strength $k  = 0.3$, the risk weight $\gamma  = 1.3$, and the mutation budget $\beta  = 0.05$. The parameter settings for the remaining competing algorithms are summarized in Table \ref{Tab3}.

\subsection{Case study}
In this subsection, the proposed model and enhanced algorithm are applied to the 2023 and 2024 O\&G exploration investment deployment cases to evaluate their effectiveness and practicality, thereby providing reproducible decision support for drilling deployment.
\begin{figure*}[htbp]
	\centering
	\includegraphics[width=\linewidth]{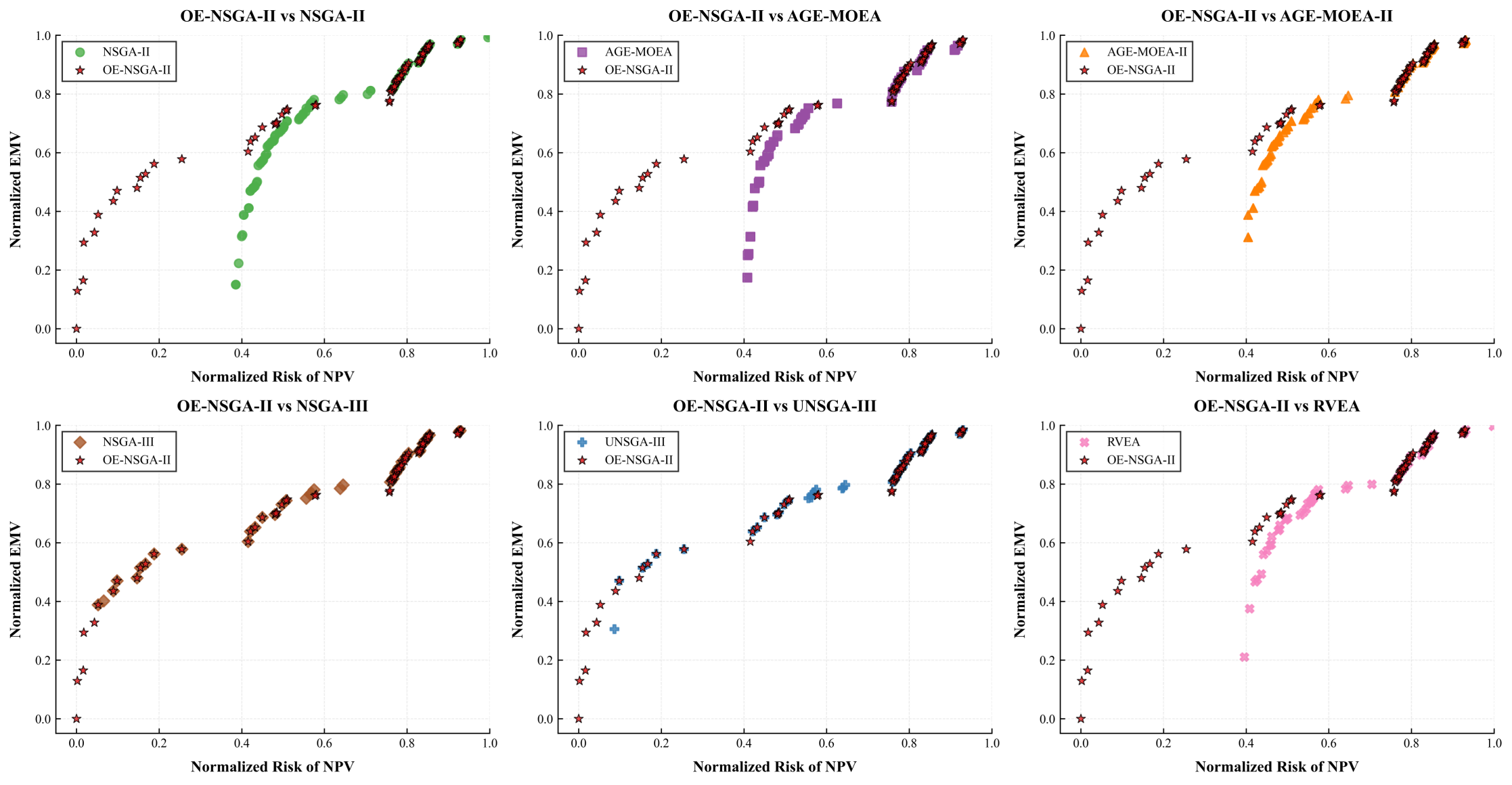} 
	\caption{Schematic diagram of the distribution of the Pareto front for each algorithm.}
	\label{fig7}
\end{figure*}
\begin{table}[htbp]
	\caption{Experimental results of the 2023 exploration investment portfolio.}\label{Tab5}
	\centering
	\setlength{\tabcolsep}{3.3pt}
	\begin{tabular}{lcclcc}
		\midrule
		Option & EMV & Risk & Option & EMV & Risk \\ 
		\hline
		Opt 1 & 5.33E+05 & 1.18E+05 & Opt 5 & 3.61E+05 & 9.08E+04 \\
		Opt 2 & 4.92E+05 & 1.12E+05 & Opt 6 & 3.62E+05 & 9.18E+04  \\
		Opt 3 & 4.85E+05 & 1.09E+05 & Opt 7 & 4.12E+05 & 9.82E+04 \\
		Opt 4 & 4.11E+05 & 9.70E+04 & Opt 8 & 5.34E+05 & 1.19E+05 \\
		\bottomrule 
	\end{tabular}
\end{table}

\subsubsection{Method validation}
In 2023, 19 trap or appraisal drilling projects must be selected from all candidates for investment, along with a certain number of reserve-providing projects as backup; in addition, some projects are designated as mandatory to support exploration in new blocks.

To better support deployment decisions while reflecting decision-makers’ preference for a return–risk trade-off, we select interpretable representative solutions using three criteria: (i) the ideal-point method to obtain a globally oriented compromise optimum; (ii) the knee-point method to identify the compromise solution with the best cost–benefit ratio (largest marginal gain); and (iii) the hypervolume-contribution method to choose the solution that maximally improves front coverage, characterizing overall performance.

Table \ref{Tab4} reports the HV, IGD, Spacing, and SC values for all algorithms. In this case, the HV of OE-NSGA-II and NSGA-III are $ 5.4877\mathrm{E}{+}09 $ and $ 5.3562\mathrm{E}{+}09 $, respectively; the IGD of OE-NSGA-II and U-NSGA-III are $ 2.6342\mathrm{E}{+}02 $ and $ 6.3544\mathrm{E}{+}02 $; and the Spacing of OE-NSGA-II and AGE-MOEA-II are $ 2.3363\mathrm{E}{+}02 $ and $ 5.9235\mathrm{E}{+}02 $, respectively. For the SC metric, the most competing algorithm is NSGA-II, with SC values of ($ 3.4694\mathrm{E}{-}01 $, $ 3.2484\mathrm{E}{-}01 $). Overall, these results indicate that OE-NSGA-II is highly competitive across all considered metrics in this case.

To compare the model-generated portfolios with the actual deployment, Table \ref{Tab5} lists eight selected solutions from the Pareto front.
\begin{figure*}[htbp]
	\centering
	\includegraphics[width=\linewidth]{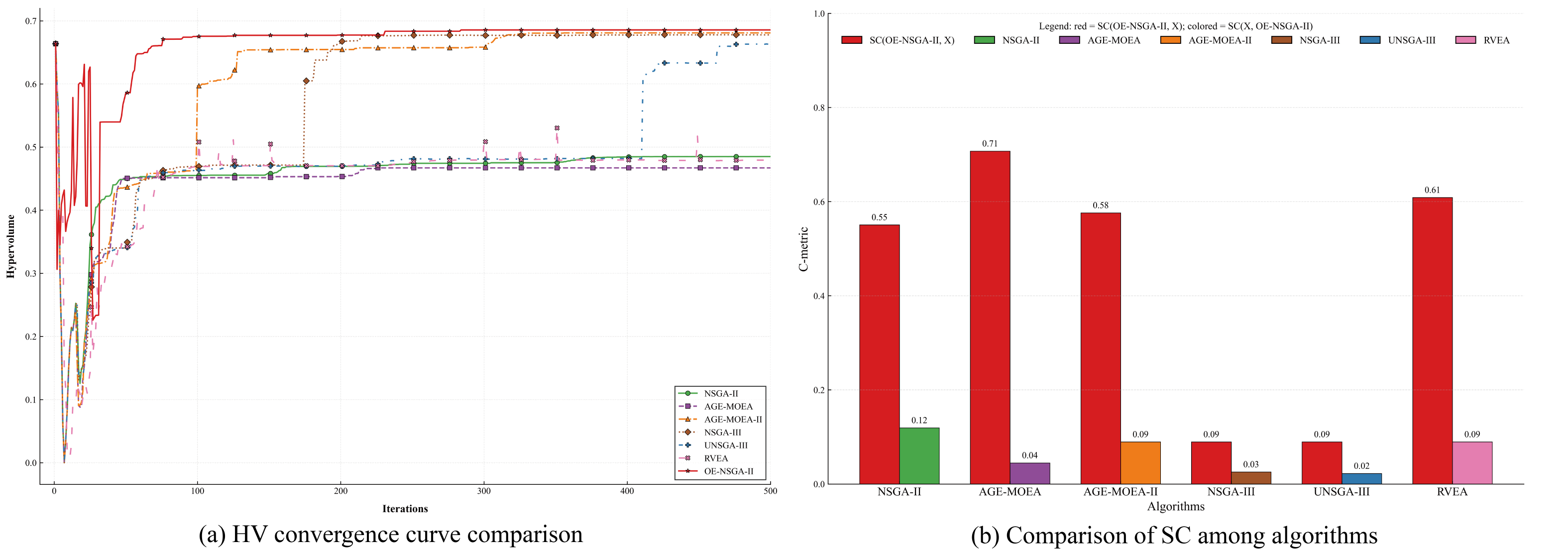} 
	\caption{Comparison of algorithms with HV and SC metrics.}
	\label{fig8}
\end{figure*}
\begin{figure*}[H!]
	\centering
	\includegraphics[width=0.8\linewidth]{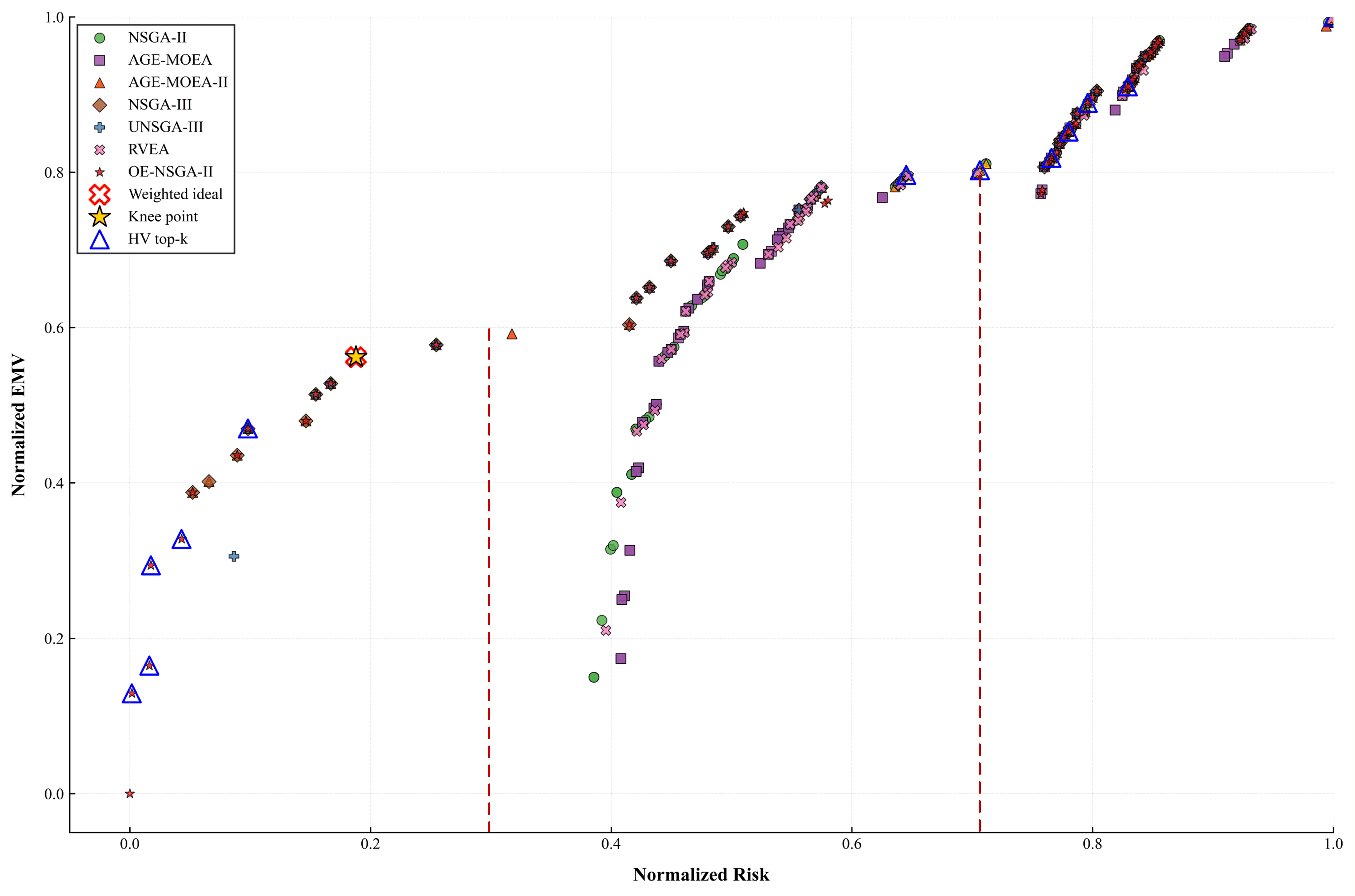} 
	\caption{Relative positions of typicale solutions on the Pareto front.}
	\label{fig9}
\end{figure*}

The risk-stratified representative solutions selected by the above procedure are shown in Fig. \ref{fig6}. As indicated in Fig. \ref{fig6} (a), the two objectives remain conflicting across risk tiers, requiring a trade-off by decision-makers. For ease of comparison, lower-risk representatives within each tier are prioritized, and the red-marked points are used as focal solutions to contrast with historical drilling selections. Fig. \ref{fig6} (b)- (c) compare eight optimized portfolios with the 2023 field deployment. Under identical budget constraints, each optimized portfolio contains no fewer projects than the 2023 plan. In terms of intersection size, the similarity between portfolios produced by our method and the 2023 deployment is approximately 58.3\%–70.8\%. While this level of agreement is attractive relative to expert-only planning, similarity alone is not a sufficient or fair evaluation. We therefore visualize the project-level differences (Fig. \ref{fig6} (c)). Discrepancies arise primarily in trap drilling projects. Beyond the overlap with the official plan, the model recommends SB14, SB8N1, SN3, XH7, YQX12, SB10, SB4, SBN82, and SB6; the plan-only set includes ADB3, QL3, SB2F1, TSA1, YQZK1, and Z1A1. Field verification shows that, apart from SB2F1 and Z1A1 (successful), the remaining projects that differ from the model’s selection failed, which supports the model’s selection quality. Notably, there are failures even within the overlap between the model and the plan—e.g., K4X1 and BST1 (designated mandatory to expand to a new area), SB8NY1, and YQX12—reflecting the inevitable risk in traps with immature geological control and high uncertainty. By contrast, for appraisal drilling and reserve-providing appraisal projects, agreement between the model and the plan is substantially higher, with no failures observed. Overall, the proposed method demonstrates effectiveness and robustness in the 2023 exploration selection.
\begin{table}[htbp]
	\caption{HV, IGD, and Spacing values obtained by all algorithms.}\label{Tab6}
	\centering
	\begin{tabular}{lccc}
		\midrule
		Algorithm & HV  & IGD  & Spacing  \\ 
		\hline
		OE-NSGA-II & 6.6982E+09  & 3.6567E+02  & 3.4323E+02  \\ 
		NSGA-III & 6.5964E+09  & 1.0656E+03  & 1.3613E+03  \\ 
		UNSGA-III & 6.4985E+09  & 9.3213E+02  & 4.2329E+03  \\ 
		NSGA-II & 5.4062E+09  & 1.1945E+03  & 1.3674E+03  \\ 
		RVEA & 5.3601E+09  & 1.6071E+03  & 4.5179E+03  \\ 
		AGE-MOEA-II & 5.3287E+09  & 1.5406E+03  & 5.2766E+02  \\ 
		AGE-MOEA & 5.2276E+09  & 1.9377E+03  & 3.1932E+03  \\ 
		\bottomrule 
	\end{tabular}
\end{table}
\begin{figure*}[htbp]
	\centering
	\includegraphics[width=\linewidth]{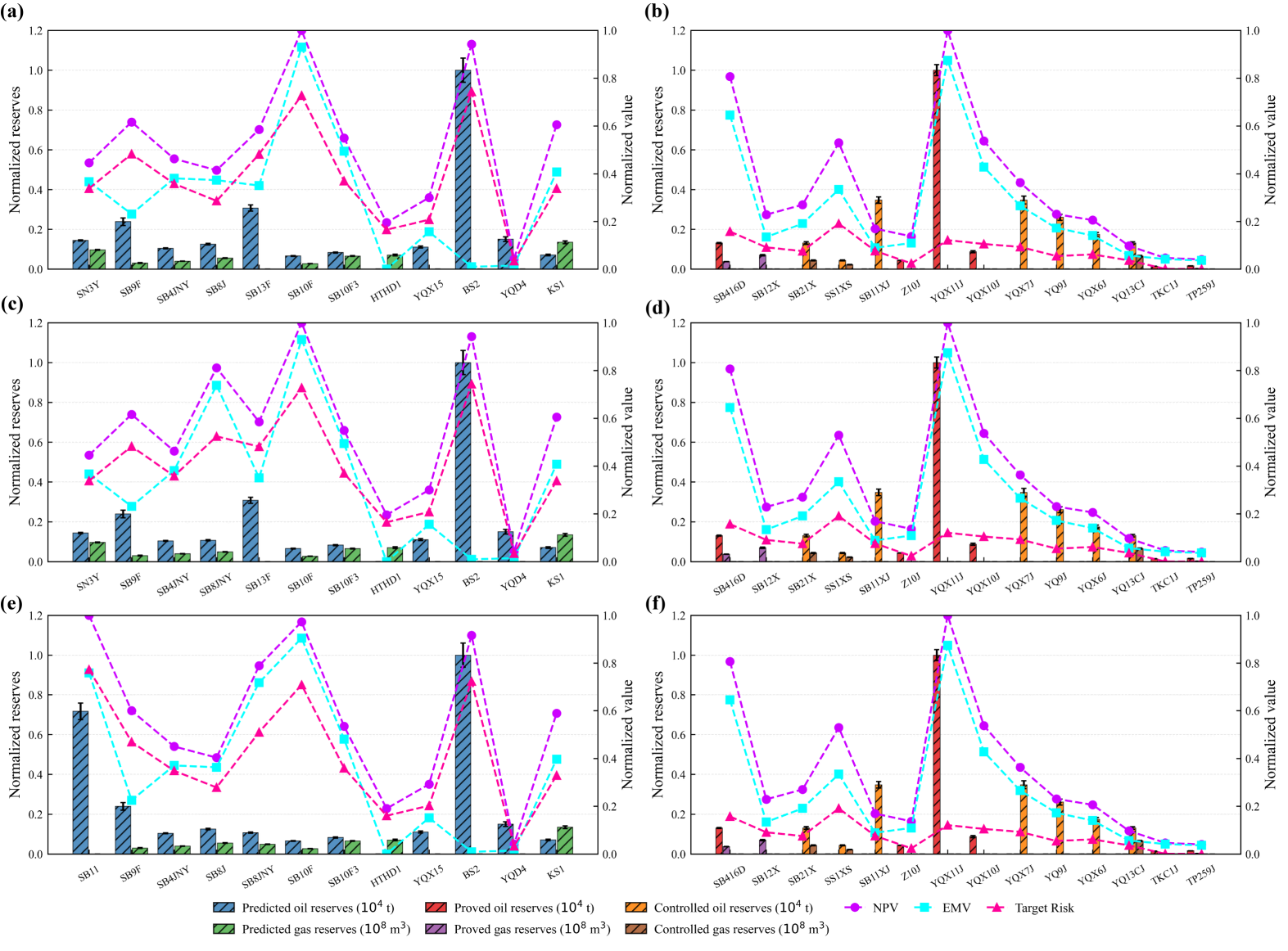} 
	\caption{Project portfolios of typical solutions by risk tier.}
	\label{fig10}
\end{figure*}

\subsubsection{Method application}
In the preceding subsection, the effectiveness of the proposed approach was validated. The model and algorithm are now applied to the 2024 exploration deployment case. All parameters and constraints are retained as before; the only difference is that the candidate project set changes in 2024, requiring optimization over a new portfolio of projects.

After applying all algorithms to the model-solving task, the resulting HV, IGD, and Spacing metrics are reported in Tables \ref{Tab6} and a comparison of the Pareto fronts is shown in Fig. \ref{fig7}. The Pareto-front comparison shows that INSGA2 achieves a superior balance between convergence and diversity, thereby offering greater decision relevance. Its front exhibits a more uniform shape with dominance advantages across multiple regions, a superiority that is particularly evident over NSGA-II, AGE-MOEA, AGE-MOEA-II, and RVEA. In the low-risk tier, OE-NSGA-II expands the feasible set (with normalized data primarily in the 0–0.25 range), offering a richer menu of options for conservative investors. In the mid-risk tier, OE-NSGA-II achieves EMV comparable to or higher than competing methods at the same risk, and lower risk at the same EMV. The solutions are also more evenly distributed, indicating stronger dominance and convergence quality, thereby better accommodating diverse risk preferences in practice. In the high-risk extreme tier (normalized return 0.7–1), OE-NSGA-II likewise approaches the theoretical bound/upper envelope, placing solutions closer to the upper hull of the efficient front and demonstrating a pronounced front-approximation advantage. Among the competing algorithms, NSGA-III and U-NSGA-III are the most competitive, each yielding desirable solution sets across the risk tiers. However, OE-NSGA-II achieves a further expansion of the low-risk region and exhibits a more uniform and more dominant Pareto set across the entire range; this advantage is particularly evident in the comparative plots.

The performance metrics in Tables \ref{Tab5} further corroborate these findings: by expanding the low-risk region, OE-NSGA-II achieves the highest HV among all competitors; its front lies closer to the upper envelope of the true Pareto front, yielding the smallest IGD; and its solutions are more uniformly distributed with clear dominance across segments, resulting in a smaller Spacing, again ranking first. Fig. \ref{fig8} presents the HV convergence curves and the SC bar chart for all algorithms. The proposed OE-NSGA-II approaches the Pareto front within approximately 80 iterations while maintaining higher front quality, whereas NSGA-III requires about 200 iterations to reach HV saturation and still underperforms the proposed method. Regarding the SC metric, OE-NSGA-II shows a pronounced advantage over NSGA-II, AGE-MOEA, AGE-MOEA-II, and RVEA, with coverage pairs of (0.55, 0.12), (0.71, 0.04), (0.58, 0.09), and (0.61, 0.09), respectively (reported as $ \text{SC} $(OE-NSGA-II, X) and $ \text{SC} $(X, OE-NSGA-II)). Although NSGA-III and U-NSGA-III remain highly competitive, OE-NSGA-II retains an overall lead. In summary, OE-NSGA-II outperforms competing algorithms in both convergence speed and solution quality.
\begin{figure*}[htbp]
	\centering
	\includegraphics[width=\linewidth]{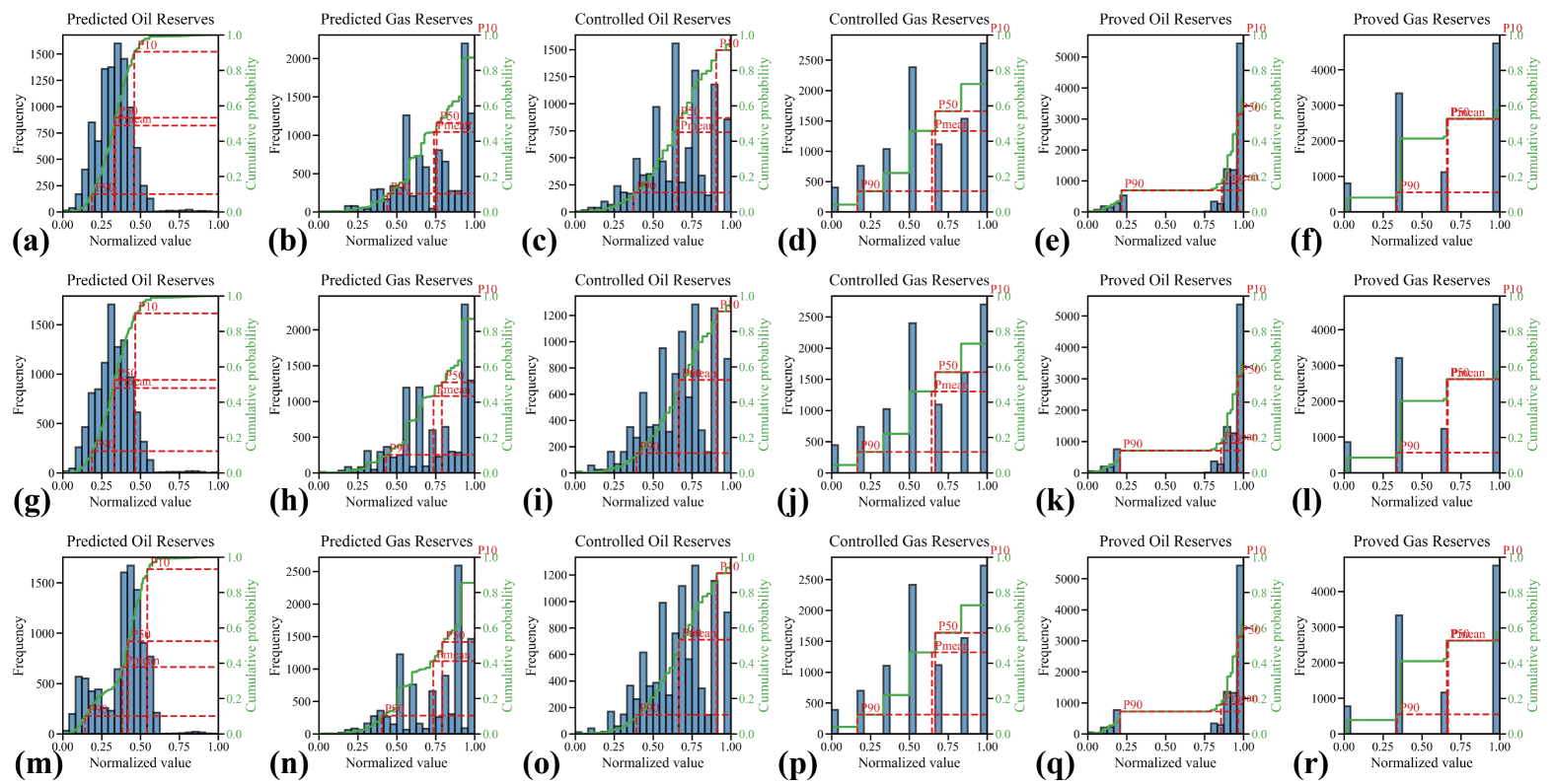} 
	\caption{Cumulative distribution functions (CDFs) of key parameters for the representative solutions.}
	\label{fig11}
\end{figure*}

\subsubsection{The typical solutions}
The Pareto front generated by the proposed optimization model yields multiple feasible solutions, whose overall distribution can be stratified into low-risk, mid-risk, and high-risk regions (see earlier discussion). The selected results obtained by the above methods are shown in Fig. \ref{fig9}, and the detailed portfolios are presented in Fig. \ref{fig10}. In addition, acknowledging uncertainty in the input parameters, their cumulative distribution functions (CDFs) are shown in Fig. \ref{fig11}, characterizing the statistical distributions of key variables and supporting subsequent robustness analysis.

Across the Pareto solutions obtained in this study, all portfolios satisfy the full set of constraints. A notable observation from the representative solutions is that selected trap projects differ appreciably across risk tiers while the overall portfolio structure remains largely stable; by contrast, appraisal projects are largely consistent across low, mid, and high risk. This indicates that portfolio risk is primarily driven by trap projects, attributable to their higher geological uncertainty and investment exposure—an outcome aligned with domain knowledge. Comparing the low-risk (Fig. \ref{fig10} (a)) and mid-risk (Fig. \ref{fig10} (c)) representatives, the key differences concentrate on SB8J and SB8JNY. Although their predicted oil/gas reserves are similar, their NPV, EMV, and target risk differ markedly. The model thus interprets these as candidates with overstated expectations, lacking stability and showing low investment reliability; even under a “high return at acceptable risk” objective, such traps materially elevate portfolio risk. A further comparison between the mid-risk (Fig. \ref{fig10} (c)) and high-risk (Fig. \ref{fig10} (e)) representatives shows that the mid-risk portfolio includes SN3Y and SB13F, whereas the high-risk portfolio includes SB11 and SB8J. SN3Y contributes to both O\&G reserves with moderate overall economic value; SB11 contributes primarily to oil reserves with relatively higher economic value but substantially higher risk, which is undesirable under stringent risk constraints. For the appraisal projects (Fig. \ref{fig10} (b), \ref{fig10} (d), \ref{fig10} (f)), geological information is complete and the economic value is pronounced. The selected candidates exhibit small overall differences in NPV and EMV, with target risk largely remaining low; consequently, this class of portfolios is both effective and robust. In summary, the risk-tiered portfolios are both engineering-plausible and interpretable: trap projects shape the portfolio’s primary risk profile, while appraisal projects provide a stable allocation backbone. Adjusting trap selections across the low–mid–high tiers enables explicit tuning of the return–risk preference.

Fig. \ref{fig11} presents the cumulative distribution functions (CDFs) of key parameters for representative solutions in the low-risk (a)-(f), mid-risk (g)-(l), and high-risk (m)-(r) tiers. The distributions of controlled oil reserves, controlled gas reserves, proved oil reserves, and proved gas reserves are broadly consistent, corroborating the preceding analysis. The main differences arise in the predicted oil and predicted gas reserves: the Pmean/P50 of predicted reserves in the low-risk tier are generally lower than those in the other tiers, whereas the mid- and high-risk portfolios tend to favor blocks with higher predicted reserves to pursue greater returns. Notably, such gains come with higher risk, which is consistent with the model’s internal mechanism.

\section{Conclusion}
We propose a multi-objective portfolio model that accounts for geological-parameter uncertainty to optimize drilling plans in O\&G exploration. The model explicitly captures the impact of geological factors under operational constraints and balances return and risk within a mean–variance framework, yielding interpretable compromise solutions for decision-makers. To address mandatory selections and combinatorial constraints, we develop an operator-enhanced NSGA-II (OE-NSGA-II) with targeted enhancements to crossover and mutation. First, a directional crossover injects the parents’ “improving direction” in objective space—approximated via dominance and objective differences—into recombination: bits aligned with the superior parent receive higher inheritance probability, while mild perturbations guide offspring to converge tangentially along the Pareto front and approach knee regions normally. This preserves diversity while improving directional search toward high-quality non-dominated solutions. Second, a structure-aware (guided) mutation stochastically prioritizes flipping high-utility bits and couples this with budget and feasibility repairs; without violating feasibility, it biases moves toward an improved return–risk trade-off, thereby strengthening the search for superior Pareto solutions. We validate the approach on a 2023 case study: the agreement between our selected portfolios and the official plan is 58.3\%–70.8\%; field checks show that model-recommended alternatives avoid failures present in the official plan, demonstrating effectiveness and practicality. We further benchmark on a 2024 case. Across HV, IGD, Spacing, SC, and the shape of the Pareto sets, OE-NSGA-II exhibits strong competitiveness relative to competing algorithms.

Future work will extend this study to mid-long term drilling deployment decisions in O\&G exploration, explicitly modeling the dynamic evolution of candidate project states as exploration progresses. Here, dynamic programming tools will be used to derive intertemporal optimal decisions under uncertainty.

\section*{Credit authorship contribution statement}
\textbf{Chao Min:} Writing – original draft, Investigation, Resources, Methodology, Formal analysis. \textbf{Junyi Cui}: Formal analysis, Software, Data curation, Writing – review \& editing, Validation. \textbf{Stanis{\l}aw Mig{\'o}rski}: Validation, Formal analysis, Visualization, Supervision. \textbf{Yonglan Xie}: Formal analysis, Validation, Supervision. \textbf{Qingxia Zhang}: Validation, Writing – review \& editing, Supervision. \textbf{Jun Peng}: Validation, Supervision.
\section*{Declaration of competing interest}

The authors declare that they have no known competing financial interests or personal relationships that could have appeared to influence the work reported in this paper.

\section*{Acknowledgments}
This work was supported by the Open Bidding for Selecting the Best Candidates Project of Southwest Petroleum University (2024CXJB11) and Sichuan Provincial Science and Technology Project (No.2025NSFTD0016).


\section*{Data availability}
The data that has been used is confidential.

\bibliographystyle{elsarticle-num}
\bibliography{ces-refs}
\end{document}